\documentclass{amsart}[12pt]
 \usepackage{geometry}        
 \usepackage{graphicx}
 \usepackage{amssymb}
 \usepackage{epstopdf}
\usepackage{amsmath}
\usepackage{amsfonts}
\usepackage{mathtools}
\usepackage{amsthm}
\usepackage{hyperref}
\usepackage{enumerate}
\usepackage{indentfirst}

\newtheorem{thm}{Theorem}
\newtheorem{prop}[thm]{Proposition}
\newtheorem{ex}[thm]{Example}
\newtheorem{lem}[thm]{Lemma}
\newtheorem{defn}[thm]{Definition}
\newtheorem{asmp}[thm]{Assumption}
\newtheorem{cor}[thm]{Corollary}
\newtheorem*{thm*}{Theorem}
\newtheorem*{thmA}{Theorem A}
\newtheorem*{thmB}{Theorem B}
\newtheorem*{thmC}{Theorem C}
\newtheorem*{thmD}{Theorem D}

\numberwithin{thm}{section}

\newcommand{\calR}{\mathcal{R}}
\newcommand{\calQ}{\mathcal{Q}}
\newcommand{\calA}{\mathcal{A}}
\newcommand{\calB}{\mathcal{B}}
\newcommand{\calL}{\mathcal{L}}

\newcommand{\calI}{\mathcal{I}}
\newcommand{\frI}{\mathfrak{I}}

\newcommand{\N}{\mathbb{N}}
\newcommand{\R}{\mathbb{R}}
\newcommand{\C}{\mathbb{C}}
\newcommand{\Z}{\mathbb{Z}}
\newcommand{\Q}{\mathbb{Q}}
\newcommand{\acl}{\operatorname{acl}}
\newcommand{\forkindep}{\mathrel{\raise0.2ex\hbox{\ooalign{\hidewidth$\vert$\hidewidth\cr\raise-0.9ex\hbox{$\smile$}}}}}
\renewcommand{\phi}{\varphi}

\newcommand{\cal}[1]{\ensuremath{\mathcal{#1}}}
\newcommand{\Cal}[1]{\ensuremath{\mathcal{#1}}}

\newcommand{\qf}{\operatorname{qf}}
\newcommand{\qftp}{\operatorname{qftp}}
\usepackage{color}

\newcommand{\ac}{\operatorname{ac}}
\newcommand{\rc}{\operatorname{rc}}
\newcommand{\dcl}{\operatorname{dcl}}

\newcommand{\tp}{\operatorname{tp}}
\newcommand{\calK}{\mathcal{K}}

\newcommand{\ind}{\operatorname{ind}}

\newcommand{\calM}{\mathcal{M}}
\newcommand{\calN}{\mathcal{N}}
\newtheorem{fact}[thm]{Fact}
\newcommand{\calO}{\mathcal{O}}
\newcommand{\calG}{\mathcal{G}}
\newcommand{\calLB}{\mathcal{L}_{\beta}}
\newcommand{\calLA}{\mathcal{L}_{\alpha}}
\newcommand{\TA}{T_{\alpha}}
\newcommand{\TB}{T_{\beta}}
\newcommand{\OR}{\overline{\mathbb{R}}}

\newcommand{\F}{\mathbb{F}}
\newcommand{\noopsort}[1]{}

\providecommand{\MR}{\relax\ifhmode\unskip\space\fi MR }

\providecommand{\href}[2]{#2}

\begin{document}

\title[Pairs of theories satisfying a Mordell-Lang condition]{Pairs of theories satisfying a Mordell-Lang condition}
\author[Block Gorman]{Alexi Block Gorman}
\author[Hieronymi]{Philipp Hieronymi}
\author[Kaplan]{Elliot Kaplan}
\email{atb2@illinois.edu}
\email{phierony@illinois.edu}
\email{eakapla2@illinois.edu}
\address{Department of Mathematics, University of Illinois at Urbana-Champaign}
\address{1409 W. Green Street, Urbana, IL 61801}
\date{\today}
\maketitle

\begin{abstract}
This paper proposes a new setup for studying pairs of structures. This new framework includes many of the previously studied classes of pairs, such as dense pairs of o-minimal structures, lovely pairs, fields with Mann groups, and $H$-structures, but also includes new ones, such as pairs consisting of a real closed field and a pseudo real closed subfield, and pairs of vector spaces with different fields of scalars. 
We use the larger generality of this framework to answer, at least in part, a couple concrete open questions raised about open cores and decidability. 
The first is, for which subfields $K \subseteq \R$ is $\R$ as an ordered $K$-vector space expanded by a predicate for $\Q$ decidable?
The second is whether there is a subfield $K$ of a real closed field that is not real closed, yet every open set definable in the expansion of the real field by $K$ is semi-algebraic.
\end{abstract}

\section{Introduction}\label{s:intro}
Pairs of structures have been widely studied in model theory, and this paper further contributes to this area. 
The goal of this paper is to describe a general framework that allows us to deduce as special cases many of the known results about pairs of structures. Among others, dense pairs of o-minimal structures as studied by van den Dries \cite{vdD98}, $H$-structures as introduced by Berenstein and Vassiliev \cite{BV16} and expansions of fields by Mann groups as discussed in van den Dries and G\"unayd\i n \cite{DG06} fall within this new framework. The larger generality of this set up will allow us to answer questions which had been outside the scope of the earlier work. Before discussing these questions and their answers, we briefly outline the new general framework.\newline

\noindent Consider a language $\calLB$ and an $\calLB$-structure $\calB$. For the moment, let $\calLA$ be a sublanguage of $\calLB$ and let $\calA$ be an $\calLA$-substructure of the $\calLA$-reduct of $\calB$. For example, consider the complex field $\C$ in the language of rings $\calL=\{0,1,+,-,\cdot\}$ and $\Gamma$ a finitely generated subgroup of $\C^{\times}$ in the sublanguage $\calL_{\frak{m}}=\{1,\cdot\}$ of multiplicative monoids. The Mordell-Lang conjecture states that for every subvariety $X$ of $\C^n$ the intersection $X\cap \Gamma^n$ is a finite union of cosets of subgroups of $\Gamma$. This implies that every $\calL$-dependence of elements in $\Gamma$ comes from a $\calL_{\frak{m}}$-dependence (see \cite[Proposition 1.1]{DG06} and Pillay \cite{Pi98}). This is the key ingredient in the proof of quantifier-elimination and model-theoretic tameness results for the pair $(\C,\Gamma)$. In this paper we will study a large class of pairs $(\calB,\calA)$ in which $\calLB$-dependence among elements in $\calA$ implies $\calLA$-dependence and, using this property, prove analogous quantifier-elimination results for $(\calB, \calA)$. Because we use this consequence of the Mordell-Lang conjecture axiomatically, we will call such pairs \textbf{ML-pairs} (for a precise definition see Section \ref{s:setup}). It is worth pointing out that we will drop the assumption that $\calLA$ is a sublanguage of $\calLB$, but this will require a more delicate definition of what a pair $(\calB,\calA)$ precisely is. We postpone this until Section \ref{s:setup}.

\subsection*{Pairs of vector spaces} Let $K$ be a subfield of $\R$. For $k \in K$, let $\lambda_k : \R \to \R$ be the function that maps $x\in \R$ to $kx$.
We denote by $\mathbb{R}_{K}$ the $K$-vector space structure on $\R$; that is the structure $\big(\R,<,+,(\lambda_k)_{k\in K}\big)$. By \cite{Hi15} the expansion of $\mathbb{R}_K$ by a predicate for $\Z$ is decidable if and only if $K$ is a quadratic field. More is true: when $K$ is not a quadratic field, then $(\mathbb{R}_K,\Z)$ defines full multiplication on $\R$ and therefore defines every open subset of $\R^n$ for every $n\in \N$. Such an expansion is as wild as can be from a model-theoretic point of view. The question was raised in \cite{Hi15} whether similar results hold when $\Z$ is replaced by $\Q$; in particular whether there is some subfield $K$ such that $(\R_K,\Q)$ is not model-theoretically well-behaved. Here we show the following.

\begin{thmA}[Corollary \ref{thma}] Every subset of $\R^n$ definable in $(\R_{K},\Q)$ is a boolean combination of sets of the form
\[
\bigcup_{\vec{q} \in \Q^m} \big\{ \vec{a} \in \R^n \ : \ (\vec{q},\vec{a}) \in X\big\},
\]
where $X\subseteq \R^{m+n}$ is definable in $\R_K$. Furthermore, every open subset of $\R^n$ definable in $(\mathbb{R}_K,\Q)$ is already definable $\mathbb{R}_K$.
\end{thmA}

\noindent Thus definable sets in $(\R_K,\Q)$ are topologically and geometrically rather tame for every subfield $K$. Furthermore, we will deduce from the proof of Theorem A that $(\mathbb{R}_K,\Q)$ is NIP (see for example Simon \cite{Si15} for a definition), and thus also exhibits strong Shelah-style model-theoretic tameness. 
Despite this model-theoretic tameness of the structure $(\mathbb{R}_K,\Q)$, its theory does not have to be decidable. For example, when $K=\R$, it is easy to see that even the theory of $\R_K$ itself is undecidable. We obtain the following characterization for when the theory of $(\mathbb{R}_K,\Q)$ is decidable.

\begin{thmB}[Theorem \ref{decidability}]
The theory of $(\mathbb{R}_K,\Q)$ is decidable if and only if 
\begin{enumerate}[(i)]
\item $K$ is a subfield of $\R$ with a computable presentation as an ordered field,
\item the question whether a finite subset of $K$ is $\Q$-linearly independent is decidable.
\end{enumerate}
\end{thmB}

\noindent Examples of such $K$ are the field of real algebraic numbers, $\Q(e^a)$ where $a\in \Q$, and $\Q(\pi)$. Note that in all of these cases, the theory of $(\R_K,\Q)$ is decidable, while the theory of $(\R_K,\Z)$ is not \cite{Le09}.

\subsection*{Pseudo real closed subfields} Let $\overline{\R}$ denote the real field. In \cite{Mi05} Miller raised the question whether for every subfield $E$ of $\R$ one of the following two statements holds:
\begin{enumerate}
\item every open set definable in $(\overline{\R},E)$ is semi-algebraic,
\item $(\overline{\R},E)$ defines $\Z$.
\end{enumerate}
As was already pointed out in \cite{Mi05}, by classical results of J. Robinson and R. Robinson, if $E$ is either a finite degree algebraic extension of $\Q$ or of the form $K(\alpha)$ with $\alpha$ transcendental over a subfield $K$, then $\Z$ is definable in just $(E,+,\cdot)$ and therefore also in $(\overline{\R},E)$. 
However, by \cite{vdD98} every open set definable in $(\overline{\R},E)$ is semi-algebraic whenever $E$ is real closed. While an answer to Miller's question is still out of reach, we are able to give the first example of subfield $E$ of $\R$ that is not real closed, but still every open set definable in $(\overline{\R},E)$ is semi-algebraic.\newline

\noindent We say that a field $K$ is \textbf{pseudo real closed} if $K$ is existentially closed in every field extension $L$ to which all orderings of $K$ extend and in which $K$ is algebraically closed.
pseudo real closed fields were first studied by Basarab \cite{Ba84} and Prestel \cite{Pr81}, and studied by van den Dries in \cite{vdD78}. Here we show the following.


\begin{thmC}[Corollary \ref{thmc}] Let $E$ be a pseudo real closed subfield of $\R$ with finitely many orderings. Then every open set definable in $(\OR,E)$ is semi-algebraic.
\end{thmC}

\noindent Since every real closed field is pseudo real closed, this generalizes the result from \cite{vdD98}. However, there are pseudo real closed subfields of $\R$ that are not real closed. Therefore Theorem C gives the desired examples of non-real closed subfields.

\subsection*{Wild theories with $P$-minimal open core} Let $\calM=(M,\dots)$ be a first order topological structure in the sense of Pillay \cite{Pi87}. The \textbf{open core of $\calM$}, denoted by $\calM^{\circ}$, is the structure $\big(M,(U)_{U \in \Cal U}\big)$, where $\Cal U$ is the collection of all open sets of all arities definable in $\calM$. Let $T^*$ be a first order topological theory in a language $\calL^*$, and let $T$ be another theory in a language $\calL$. We say that \textbf{$T$ is an open core of $T^*$} if for every $\calN\models T^*$ there is $\calM \models T$ such that $\calN^{\circ}$ is interdefinable with $\calM$. The notion of an open core of a theory was introduced in Dolich, Miller, and Steinhorn \cite{DMS10} for theories extending the theory of dense linear orders, generalizing earlier work of Miller and Speissegger on expansions of the real line \cite{MS99}.\newline

\noindent Hieronymi, Nell and Walsberg \cite{HNW17} investigated the question of whether there are any tameness conditions that can be imposed on the open core (such as o-minimality) such that the whole theory satisfies some (possibly weaker) form of the model-theoretic tameness. In that paper a rather strong negative answer was given in case the open core is o-minimal. However, the same question for theories with \textbf{$P$-minimal} open core was left open. The notion of $P$-minimality was introduced in Haskell and Macpherson \cite{HM97}, where it was developed as an analog to o-minimality for $p$-adically closed fields. Here we use our general framework to answer the above tameness question for $P$-minimal open cores.

\begin{thmD}[Section \ref{s:padic}] Let $T$ be the theory of the $p$-adic field $\Q_p$, and let $T'$ be a consistent theory. Then there is a complete theory $T^*$ extending $T$ such that
\begin{enumerate}
\item $T^*$ interprets a model of $T'$,
\item $T$ is an open core of $T^*$,
\end{enumerate}
\end{thmD}

\noindent Since $\Q_p$ is $P$-minimal, this result rules out that the property of having an $P$-minimal open core has any consequences in terms of model-theoretic tameness of the whole theory.

\subsection*{Acknowledgements}
The authors would like to thank Lou van den Dries, Erik Walsberg, Allen Gehret, and Minh Tran for their thoughts and conversations related to this paper.
The first author was partially supported by a DOE GAANN fellowship. The second author was partially supported by NSF grant DMS-1654725. 

\subsection*{Notation and conventions} We will use $m,n$ for natural numbers and $\kappa$ for a cardinal. Let $X,Y$ be sets. We denote the cardinality of $X$ by $|X|$. If $Z\subseteq X \times Y$ and $x\in X$, then $Z_x$ denotes the set $\big\{ y\in Y \ : \ (x,y) \in Z\big\}$. If $\vec{z}=(z_1,\dots,z_n)$, we sometimes write $X\vec{z}$ for $X\cup \{z_1,\dots,z_n\}$, and $XY$ for $X\cup Y$.\newline

\noindent Let $\calL$ be a language and $T$ an $\calL$-theory. We set $|T|:=|\calL|$ if $\calL$ is infinite and $|T|:= \aleph_0$ otherwise. Let $\cal M \models T$ and $C\subseteq M$. As a matter of convenience, we view constant symbols in $\calL$ as nullary functions. We use $\calL$-definable to mean $\calL$-definable without parameters and we use $\calL(C)$-definable (or $\calL$-definable over $C$) to indicate $\calL$-definability with parameters from $C$. The same conventions hold for $\calL$-formulas and $\calL$-types. For an $\calL(M)$-formula $\phi(\vec{x})$, we write $\phi(\calM)$ to denote the $\calL(M)$-definable subset of $M^{|\vec{x}|}$ defined by $\phi$.\newline

\noindent Let $b\in M^n$. Then we write $\tp_{\calL}(b|C)$ for the $\calL$-type of $b$ over $C$ computed in $\calM$. Types are always assumed to be complete and realizable. Let $p$ be an $\calL(C)$-type. We let $\qf (p)$ denote the set of quantifier-free formulas in $p$.
Let $\calN$ be another model of $T$ and $D\subseteq N$. If $\iota: C \to D$ is a partial $\calL$-isomorphism, then we denote by $\iota p$ the set of formulas $\phi\big(\vec{x},\iota(\vec{c})\big)$ such that $\phi(\vec{x},\vec{c}) \in p$. 
This is indeed a type (i.e. it is realizable) if $\iota$ is an $\calL$-elementary map.\newline


\section{Setup}\label{s:setup}

Consider a language $\calLB$ and a consistent $\calLB$-theory $\TB$. Let $\calLA$ be another language whose function symbols are all in $\calLB$. Let $\TA$ be a consistent $\calLA$-theory. We denote the intersection of $\calLA$ and $\calLB$ by $\calL$. Let $\calL^2 = \calLB \cup \calLA \cup \{ A \}$ where $A$ is a unary predicate symbol not contained in $\calLB \cup \calLA$.\newline

\noindent Let $\theta$ be an $\calLA$-formula. We define the $\calLA\cup\{A\}$-formula $\theta_A$ by relativizing all of the quantifiers in $\theta$ to the predicate $A$. More precisely, we define $\theta_A$ recursively as follows:\\
\hspace{1mm} \hspace*{20mm} $\theta_A := t_1(\vec{x})=t_2(\vec{x}), \text{ if } \theta \text{ is } t_1(\vec{x})=t_2(\vec{x})$ and $t_1,t_2$ are $\calLA$-terms \\
\hspace{1mm} \hspace*{20mm} $ \theta_A :=  R\big(t_1( \vec{x}),\ldots,t_n(\vec{x})\big), \text { if } \theta \text{ is } R\big( t_1(\vec{x}),\ldots,t_n(\vec{x}) \big)\text{ where $t_1,\ldots,t_n$ are $\calLA$-terms }\\
\hspace{1mm} \hspace*{31mm}\text{and $R$ is a relation symbol in $\calLA$}$ \\
\hspace{1mm} \hspace*{20mm} $ \theta_A := \neg \theta _A ', \text{ if } \theta \text{ is } \neg \theta ' $ \\
\hspace{1mm} \hspace*{20mm} $ \theta_A := \theta'_A \land \theta ''_A, \text{ if } \theta \text{ is } \theta' \land \theta '' $ \\
\hspace{1mm} \hspace*{20mm} $ \theta_A := \exists x \big(A(x) \land \theta'_A\big), \text{ if } \theta \text{ is } \exists x \theta ' $\\
\hspace{1mm} \hspace*{20mm} $ \theta_A := \forall x \big(A(x) \rightarrow \theta'_A\big), \text{ if } \theta \text{ is } \forall x \theta ' $.\\
Set
\[
A(\TA) := \{ \phi_A \ : \ \TA \models \phi \}.
\]
We denote by $T^2$ the $\calL^2$-theory extending $\TB \cup A(\TA)$ by the following schemas of $\calL^2$-sentences:
\begin{enumerate}
\item[(T1)] for each function symbol $f\in \calLA$ of arity $n$
\[
\forall x_1 \dots \forall x_n \Big(\big(\bigwedge_{i=1}^n A(x_i)\big) \rightarrow A\big(f(x_1,\dots,x_n)\big)\Big),
\]
\item[(T2)] for each relation symbol $R\in \calLA\setminus \calLB$ of arity $n$,
\[
\forall x_1 \dots \forall x_n \Big(R(x_1,\dots,x_n) \rightarrow \bigwedge_{i=1}^n A(x_i)\Big).
\]
\end{enumerate}
\noindent Suppose that $T^2$ has a model $\calM$. We denote the reduct of $\calM$ to $\calLB$ by $\calB_{\calM}$. Set $A_{\calM} := \big\{x \in M:\calM\models A(x)\big\}$. By (T1) we have that $A_{\calM}$ is an $\calLA$-substructure of the reduct of $\calM$ to $\calLA$. We denote this substructure by $\calA_{\calM}$.
We remark that for $T^2$ to be consistent, it is necessary that for every $\calL$-sentence $\phi$ it holds that $T_{\beta,\forall} \vdash \phi \implies \TA \vdash \phi$, and it is sufficient for both $\TB$ and $\TA$ to imply precisely the same $\calL$-sentences.
However, these do not completely characterize when $T^2$ is consistent, as we will see in the examples sections.

\begin{lem}\label{lem:uniquem} Let $\calM,\calM'\models T^2$. If $\calB_{\calM} = \calB_{\calM'}$ and $\calA_{\calM} = \calA_{\calM'}$, then $\calM = \calM'$.
\begin{proof}
Since $\calB_{\calM} = \calB_{\calM'}$, the two models $\calM$ and $\calM'$ have the same underlying set $M$. It is left to show that every symbol in $\calL^2$ is interpreted the same way in $\calM$ and $\calM'$. 
It is immediate that all symbols in $\calLB$ are interpreted equally. Furthermore, $A_{\calM}=A_{\calM'}$, because these are the underlying sets of $\calA_{\calM}$ and $\calA_{\calM'}$, and $\calA_{\calM} = \calA_{\calM'}$. 
It remains to consider symbols in $\calLA\setminus \calLB$. 
Since every function symbol in $\calLA$ is also in $\calLB$, we can reduce to relation symbols (recall that we view constants as nullary functions). 
Let $R$ be a relation symbol in $\calLA\setminus \calLB$. By (T2), $R_{\calM} = R_{\calA}$ and $R_{\calM'} = R_{\calA'}$. Thus $R_{\calM} = R_{\calM'}$ because $\calA_{\calM} = \calA_{\calM'}$.
We note that the assumption that every function symbol in $\calLA$ is also in $\calLB$ is necessary: if $\calLA$ contains a function symbol $f$ which is not in $\calLB$ then one can come up with examples where $f_{\calM}$ and $f_{\calM'}$ disagree at some tuple not in $A$.
\end{proof}
\end{lem}

\noindent From now on, when we write $(\calB,\calA) \models T^2$, we mean that there is a model $\calM$ of $T^2$ such that $\calB=\calB_{\calM}$ and $\calA=\calA_{\calM}$. When we refer to the pair $(\calB,\calA)$, we are referring to this model $\calM$. By Lemma \ref{lem:uniquem} the two structures $\calB$ and $\calA$ determine $\calM$ uniquely. We let $B$ denote the underlying set of $\calB$ and we let $A$ denote the underlying set of $\calA$.

\begin{lem}
\label{relisequiv}
Let $(\calB,\calA)\models T^2$. Let $\vec{a} \in A^n$ and let $\phi$ be an $\calLA$-formula. Then
\[
(\calB,\calA) \models \phi_A(\vec{a}) \hbox{ if and only if } \calA \models \phi(\vec{a}).
\]
In particular, $\calA\models \TA$.
\begin{proof}
This follows by a straightforward induction on $\calLA$-formulas.
\end{proof}
\end{lem}

\noindent For a tuple $\vec{a}$ from $A$ and $C \subseteq A$, we use $\tp_{\calLA} (\vec{a} |C)$ to denote the collection of all $\calLA(C)$-formulas $\psi(\vec{x})$ such that $\calA \models \psi(\vec{a})$.
Given an $\calLA(C)$-type $p(\vec{x})$, we let $p_A(\vec{x}) = \big\{\psi_A(\vec{x}):\psi \in p\big\}$.
We observe by the above lemma that $\calA \models p( \vec{a} )$ if and only if $(\calB, \calA) \models p_A (\vec{a})$.
Note also that if $\phi$ is a quantifier-free $\calLA$-formula, then $(\calB,\calA) \models \phi_A(\vec{a})$ if and only if $(\calB,\calA) \models \phi(\vec{a})$. We will use this fact often.


\subsection{ML-pairs}
From now on we assume that $\TB$ is \textbf{geometric}; that is $\TB$ eliminates the $\exists^{\infty}$ quantifier and the algebraic closure operator $\acl$ defines a pregeometry in every model of $\TB$. Let $\calB$ be a model of $\TB$. Let $X,Y,Z$ be subsets of $B$. We say that $X$ and $Y$ are \textbf{independent over $Z$} -- written as $X \forkindep_{Z} Y$ -- if every subset of $X$ that is $\acl$-independent over $Z$ is also $\acl$-independent over $YZ$. The following lemma is often useful:

\begin{lem}
\label{finitefree}
Let $X,Y \subseteq B$ and suppose that $X \forkindep_{X \cap Y} Y$. Then for every finite $X_0\subseteq X$ there is a finite $X_1 \subseteq X$ such that $X_0 \subseteq X_1$ and $X_1 \forkindep_{X_1 \cap Y} Y$.
\begin{proof}
There are only finitely many subsets of $X_0$ that are $\acl$-independent over $X_0\cap Y$, and each of these subsets has empty intersection with $Y$. If a subset $Z\subseteq X_0$ is $\acl$-independent over $X_0 \cap Y$ but not over $Y$, then there is a finite subset $Y_Z\subseteq Y$ containing $X_0 \cap Y$ such that $Z$ is $\acl$-dependent over $Y_Z$. As $X \forkindep_{X\cap Y}Y$, we can assume that $Y_Z \subseteq X \cap Y$. Let
\[
X_1:= X_0 \cup \bigcup\big\{Y_Z:Z \subseteq X_0\text{ is $\acl$-independent over $X_0\cap Y$ but not over }Y\big\}.
\]
Then $X_1$ is finite and $X_0 \subseteq X_1\subseteq X$. The reader can easily check that $X_1\forkindep_{X_1\cap Y}Y$.
\end{proof}
\end{lem}

\noindent For a model $(\calB,\calA)\models T^2$, we use $\acl$ to denote the algebraic closure in $\calB$ and we use $\dcl$ to denote the definable closure in $\calB$. 
For a tuple $\vec{c} \in B^n$, we let $\vec{c}_\alpha$ be the subtuple of $\vec{c}$ consisting of the components of $\vec{c}$ belonging to $A$. If $\vec{x}$ is a tuple of variables, we let $\vec{x}_\alpha$ be a subtuple of variables (which may be empty or equal to $\vec{x}$). We think of $\vec{x}_\alpha$ as the part of $\vec{x}$ which ranges over $A$.

\begin{defn}\label{defn:MLT}
Let $T \supseteq T^2$ be an $\calL^2$-theory. A \textbf{Mordell-Lang challenge (for $T$)} is a tuple
\[
\big(p(\vec{x}_\alpha),q(\vec{x}),\phi(\vec{x},y),\psi(\vec{x}_\alpha,y)\big)
\]
such that 
\begin{itemize}
\item $p$ is a complete $\calLA$-type and $q$ is a complete $\calLB$-type,
\item $\phi$ is an $\calLB$-formula and $\psi$ is an $\calLA$-formula,
\item $q(\vec{x}) \models \exists^{<\infty}y\phi(\vec{x},y)$.
\end{itemize}
A \textbf{contender} to a Mordell-Lang challenge is a tuple $\big((\calB,\calA),\vec{c}\ \big)$ where $(\calB,\calA)\models T$ and where $\vec{c}$ is a tuple in $B$ such that $\vec{c}$ realizes $q$, the subtuple $\vec{c}_\alpha$ realizes $p_A$, and $\vec{c}\forkindep_{\vec{c}_\alpha}A$.
A \textbf{solution} to a Mordell-Lang challenge is a tuple $\big((\calB,\calA),\vec{c},a\big)$ such that $\big((\calB,\calA), \vec{c}\ \big)$ is a contender, $a \in A$, and 
\[
(\calB,\calA) \models  \phi(\vec{c},a)\wedge \psi_A(\vec{c}_\alpha,a).
\]
A Mordell-Lang challenge is \textbf{solvable} if it has a solution.
\end{defn}


\begin{defn}\label{defn:AE}
Let $T \supseteq T^2$ be an $\calL^2$-theory. We say that $T$ satisfies the \textbf{Mordell-Lang condition} if for every solvable Mordell-Lang challenge $(p,q,\phi,\psi)$ for $T$ and for every contender $\big((\calB,\calA), \vec{c}\ \big)$, there is $a \in A$ such that $\big((\calB,\calA),\vec{c},a\big)$ is a solution.
\end{defn}

\noindent
We are now ready to define a Mordell-Lang theory of pairs.

\begin{defn}\label{defn:ML}
An $\calL^2$-theory $T$ is a \textbf{Mordell-Lang theory of pairs} (or short: \textbf{ML-theory}) if
\begin{enumerate}
\item $T$ extends $T^2$,
\item $T$ satisfies the Mordell-Lang condition,
\item for every $\kappa$-saturated model $(\calB,\calA) \models T$ where $\kappa >|T^2|$, every $C \subseteq B$ with $|C|<\kappa$, and every non-algebraic unary $\calLB(C)$-type $q(x)$ the following conditions hold:
\begin{enumerate}
\item (Density)
if $p(x)$ is a unary $\calLA(A\cap C)$-type such that $q \models \qf (p|_{\calL})$, where $p|_{\calL}$ restricts to only $\calL(A\cup C)$-formulas, then there is $a \in A$ realizing $p_A \cup q$.
\item (Codensity)
there is $b \in B\setminus \acl(A\cup C)$ realizing $q$.
\end{enumerate}
\end{enumerate}
A model $(\calB,\calA)$ of an ML-theory is called an \textbf{ML-pair}.
\end{defn}

\noindent The density and codensity conditions are inspired by the extension and coheir properties used by Berenstein and Vasseliev \cite{BV10} to axiomatize lovely pairs of geometric theories. In the case that $\calB$ has a definable topology, these don't correspond exactly to density and codensity of $A$ in $B$, but they are related. We will present examples of ML-theories in the next subsection.

\subsection{Known examples}
Here we describe three well-known classes of theories which fit into our framework. In Sections \ref{s:oml}-\ref{s:padic} we will present three classes of structures that have not been studied before, but also fall within this new setup.

\subsubsection*{Lovely pairs} Let $T_\beta$ be a geometric theory with quantifier elimination in the language $\calLB$ and set $\calLA := \calLB$ and $T_\alpha := T_\beta$. Let $T_P\supseteq T^2$ be an $\calL^2$ theory such that $T_P$ satisfies the density and codensity conditions in Definition \ref{defn:ML} and such that for any $(\calB,\calA) \models T_P$, the set $A$ is algebraically closed in $\calB$. Then $\calA$ is an elementary substructure of $\calB$ in every model of $T_P$ and any $|T^2|^+$-saturated model of $T_P$ is a \textbf{lovely pair of models of $T_\beta$}. These lovely pairs are axiomatized in Theorem 2.10 in \cite{BV10}, and their theory is studied extensively in the same paper.

\begin{prop}\label{lovelypairs}
The theory $T_P$ is an ML-theory.
\begin{proof}
By definition, $T_P$ satisfies conditions (1) and (3) in Definition \ref{defn:ML}. It remains to check that $T_P$ satisfies the Mordell-Lang condition. Let $(p,q,\phi,\psi)$ be a Mordell-Lang challenge, suppose that $\big((\calB,\calA),\vec{c},a\big)$ is a solution, and let $\big((\calB',\calA'), \vec{d}\ \big)$ be a contender. Then since $\phi(\vec{c},y)$ is an algebraic formula, we have that $a \in \acl(\vec{c})\cap A$. Using the fact that $\vec{c}\forkindep_{\vec{c}_\alpha}A$, we have $a \in \acl(\vec{c}_\alpha)\cap A$. Since $\calLB = \calLA$, we may assume that $\psi(\vec{c}_\alpha,y)$ is algebraic (if not, then replace $\psi$ with an algebraic formula which implies $\psi$). We may also assume that $\psi$ is quantifier-free, so $\psi = \psi_A$. We have $\calA \models \psi(\vec{c}_\alpha,a)$, so since $\calA $ is an $\calLB$-elementary substructure of $\calB$, we also have that $\calB \models \psi(\vec{c}_\alpha,a)$. Therefore,
\[
q(\vec{x})\models \exists y \big(\phi(\vec{x},y) \wedge \psi(\vec{x}_\alpha,y)\big).
\] 
Since $\vec{d}$ realizes $q$, there is some $a' \in B'$ such that 
$\calB' \models \phi(\vec{d},a') \wedge \psi(\vec{d}_\alpha,a').$
Since $A'$ is algebraically closed and $\psi$ is algebraic, we must have $a' \in A'$ and so $\big((\calB',\calA'),\vec{d},a'\big)$ is a solution.
\end{proof}
\end{prop}

\subsubsection*{Expansions by $\acl$-independent sets}
Let $T_\beta$ be a geometric theory and let $T_\alpha$ extend the theory of an infinite set. In particular note that we do not require $\calLA$ to be empty. Let $T_{\ind}\supseteq T^2$ be an $\calL^2$-theory that satisfies condition (3) in Definition \ref{defn:ML} and includes the sentence
\[
 \forall x_1 \dots x_n \Big(\big(\bigwedge_{i=1}^n A(x_i) \wedge \exists^{<\infty}y \ \phi(x_1,\dots,x_n,y)\big)\rightarrow \forall y \big(\bigwedge_{i=1}^n (x_i \neq y) \wedge A(y) \rightarrow \neg\phi(x_1,\dots,x_n,y)\big)\Big)
\]
for each $n$ and each $\calLB$-formula $\phi(x_1,\dots,x_n,y)$. This last axiom implies that $A$ is an $\acl$-independent set in every model $(\calB,\calA)$ of $T_P$. Furthermore, whenever $\calLA=\emptyset$, it follows easily that every $|T^2|^+$-saturated model of $T_{\ind}$ is an \textbf{$H$-structure}, as defined in \cite{BV16}.

\begin{prop}
\label{ind} The theory $T_{\ind}$ is an ML-theory.
\begin{proof}
By assumption, $T_{\ind}$ satisfies conditions (1) and (3) in Definition \ref{defn:ML}. It remains to check that $T_{\ind}$ satisfies the Mordell-Lang condition. Let $(p,q,\phi,\psi)$ be a Mordell-Lang challenge, suppose that $\big((\calB,\calA),\vec{c},a\big)$ is a solution, and let $\big((\calB',\calA'), \vec{d}\ \big)$ be a contender. Then since $a \in \acl(\vec{c})\cap A$ and $\vec{c}\forkindep_{\vec{c}_\alpha}A$, we have $a \in \acl(\vec{c}_\alpha)\cap A$. Since $A$ is $\acl$-independent, it must be the case that $a$ is a component of $\vec{c}_\alpha$. Letting $a'$ be the corresponding component of $\vec{d}_\alpha$, we have that $(\calB',\calA') \models \phi(\vec{d},a')\wedge \psi_A(\vec{d}_\alpha,a')$.
\end{proof}
\end{prop}

\subsubsection*{Algebraically closed fields with a Mann subgroup}
Let $L$ be a field and let $\Gamma$ be an infinite multiplicative subgroup of $L^\times$. We denote the prime field of $L$ by $\F$.
\begin{defn}\label{def:mann}
We say that $\Gamma$ has the \textbf{Mann property} if for every $\vec{q} = (q_1,\ldots,q_n)\in (\F^{\times})^n$ there are only finitely many tuples $\vec{\gamma}=(\gamma_1,\ldots,\gamma_n) \in \Gamma^n$ such that $\sum_{i = 1}^n q_i\gamma_i = 1$ and $\sum_{i \in I} q_i\gamma_i \neq 0$ for every nonempty $I \subseteq \{1,\ldots,n\}$. Such a tuple $\vec{\gamma}$ is called a \textbf{non-degenerate solution} to the $\F$-linear equation $\sum_{i = 1}^n q_ix_i = 1$.
\end{defn}

\noindent Many interesting multiplicative subgroups of fields have the Mann property. For instance, if $\Gamma$ has finite rank and $L$ is of characteristic 0, then $\Gamma$ has the Mann property. Pairs of fields with Mann subgroups are studied extensively in \cite{DG06}, and in the following we show that this work fits under the umbrella of ML-theories. \newline

\noindent From now on we assume that $L$ is algebraically closed and $\Gamma$ is a subgroup of $L^\times$ with the Mann property with $[\Gamma:\Gamma^n]<\infty$ for each $n \geq 1$. We will consider the case where $L$ is real-closed and $\Gamma$ is divisible in Section \ref{s:oml}. We axiomatize the pair $(L,\Gamma)$ as follows: set $\calLA:= \big\{1, \cdot, x\mapsto x^{-1}, (\gamma)_{\gamma \in \Gamma} \big\}$ and consider $\Gamma$ as an $\calLA$-structure in the natural way. Let $\TA$ be the $\calLA$-theory of $\Gamma$. Set $\calLB := \big\{0,1, \cdot,+,- , x\mapsto x^{-1}, (\gamma)_{\gamma \in \Gamma}\big\}$ and let $\TB$ be the $\calLB$-theory of $L$ (with $0^{-1}:= 0$). We let $T^{\ac}_\Gamma \supseteq T^2$ be the theory stating that for $(\calK,\calG)\models T^{\ac}_\Gamma$ and for every $\F$-linear equation $\sum_{i = 1}^n q_ix_i = 1$, each non-degenerate solution in $\calG$ is one of the solutions in $\Gamma$. Since there are only finitely many non-degenerate solutions in $\Gamma$, such an $\calL^2$-theory is indeed axiomatizable, as observed in \cite{DG06}.\newline

\noindent In order to show that $T^{\ac}_{\Gamma}$ satisfies the Mordell-Lang condition, we rely on the following Lemma, which is an immediate corollary of the proof of \cite[Proposition 5.8]{DG06}.
\begin{lem}\label{lem:MLforMann}
Let $\phi(\vec{x})$ be an $\calLB$-formula. Then there is an $\calLA$-formula $\psi(\vec{x})$ such that
\[
(\calK,\calG)\models \phi(\vec{a}) \text{ if and only if }\calG\models  \psi\big(\vec{a})
\]
for all all models $(\calK,\calG) \models T^{\ac}_{\Gamma}$ and all tuples $\vec{a}$ from $G$.
\end{lem}

\noindent For $(\calK,\calG)\models T^{\ac}_\Gamma$ and $C \subseteq K$, we let $\F(C)$ denote the subfield of $\calK$ generated by $C$. For a subgroup $E$ of $\calG$, we say that $E$ is \textbf{$\calLB$-existentially closed in $\calG$} if for each quantifier-free $\calLB(E)$-formula $\phi(\vec{x})$, if there is $\vec{a} \in G^{|\vec{x}|}$ such that $\F (G) \models \phi (\vec{a})$, then there is some $\vec{e} \in E^{|\vec{x}|}$ such that $\F (E) \models \phi (\vec{e})$. By Lemma 3.3 in \cite{DG06}, if $E$ is $\calLB$-existentially closed in $\calG$ then $\F(G)$ is a regular extension of $\F(E)$.

\begin{lem}\label{Mann}
The theory $T^{\ac}_\Gamma$ satisfies the Mordell-Lang condition.
\begin{proof}
Let $(p,q,\phi,\psi)$ be a Mordell-Lang challenge, suppose that $\big((\calK,\calG),\vec{c},a\big)$ is a solution, and let $\big((\calK',\calG'), \vec{d}\ \big)$ be a contender.
Let $\iota$ be the function mapping $\vec{c}$ to $\vec{d}$ componentwise, so $\iota$ maps $\vec{c}_\alpha$ to $\vec{d}_\alpha$, and let $\iota'$ denote the restriction of $\iota$ to $\vec{c}_\alpha$. Then $\iota$ is $\calLB$-elementary and $\iota'$ is $\calLA$-elementary. Take an $\calLB$-existentially closed subgroup $E$ of $\calG$ containing $\vec{c}_\alpha$ and $a$ and extend $\iota'$ to a $\calLA$-elementary map $\tilde{\iota}:E \to G'$. Set $E' := \tilde{\iota}(E)$ and set $a' := \tilde{\iota}(a)$. Then $\calG' \models \psi(\vec{d}_\alpha,a')$ and, by Lemma \ref{lem:MLforMann}, $E'$ is $\calLB$-existentially closed in $\calG'$. It remains to show that $\calK' \models \phi(\vec{d},a')$.  

Since $E \subseteq G$, and $\vec{c} \forkindep_{\vec{c}_\alpha} G$, we have that $E\vec{c} \forkindep_{E} G$ and so $\F(E\vec{c}) \forkindep_{\F (E)} \F (G)$. Since $E$ is $\calLB$-existentially closed in $\calG$, we also have that $\F(G)$ is a regular extension of $\F(E)$ and we conclude by \cite[p. 367]{La95} that $\F( E\vec{c})$ and $\F(G)$ are linearly disjoint over $\F (E)$. Likewise, $\F( E'\vec{d})$ and $\F(G')$ are linearly disjoint over $\F (E')$. Thus, there is an $\calLB$-isomorphism $\tilde{\iota}':\F (E\vec{c}) \xrightarrow{\sim} \F (E'\vec{d})$ which extends both $\iota$ and $\tilde{\iota}$. As $\TB$ admits quantifier elimination in the language $\calLB$, we may assume that $\phi$ is quantifier-free and so $\calK' \models \phi(\vec{d},a')$.  
\end{proof}
\end{lem}

\begin{lem}\label{Mannalg}
Let $\calG \models \TA$, let $\vec{c}$ be a tuple from $G$, and let $\phi(x)$ be an $\calLA(\vec{c})$-formula such that $\phi(\calG)$ is finite. Then there is a quantifier-free $\calLA(\vec{c})$-formula $\psi(x)$ such that  $\psi (\calG)$ is a finite set containing $\phi(\calG)$.
\begin{proof}
We consider the expansion of $\calG$ by predicates $D_n$ where 
$$D_n(\calG) = \{a \in G: h^n = a \text{ for some } h \in G\}.$$
By Szmielew's quantifier elimination for abelian groups \cite{Sz55}, $\calG$ admits quantifier elimination in this language. 
We may assume that $\phi(x)$ is equivalent to a disjunction of formulas of the form
$$\psi(x) \wedge D_{n_1}\big(x^{m_1}t_1(\vec{c})\big) \wedge \ldots \wedge D_{n_k}\big(x^{m_k}t_k(\vec{c})\big)$$
where  $\psi$ is a quantifier-free $\calLA$-formula, $m_i,n_i$ are natural numbers, and $t_i$ is an $\calLA$-term for each $i$. Assume that $\phi$ is equivalent to just one disjunct of this form. 
By letting $n$ be the least common multiple of $n_1,\ldots,n_k$ and raising $x^{m_i}t_i(\vec{c})$ to the power $n/n_i$, we may further assume that $n_1,\ldots,n_k$ are all the same. 
We note that if $D_{n}\big(x^{m_i}t_i(\vec{c})\big)$ holds for some $x \in G$, then $D_{n}\big(y^{m_i}t_i(\vec{c})\big)$ holds for all $y \in \calG^nx$ and that $\calG^n x$ is infinite (as $[\calG:\calG^n]$ is finite). 
Thus, 
$$D_{n_1}\big(x^{m_1}t_1(\vec{c})\big) \wedge \ldots \wedge D_{n_k}\big(x^{m_k}t_k(\vec{c})\big)$$
must define an infinite set and so $\psi(\calG)$ must be a finite set containing $\phi(\calG)$.
\end{proof}
\end{lem}

\begin{prop}
The theory $T^{\ac}_\Gamma$ is an ML-theory.
\begin{proof}
By Lemma \ref{Mann} the theory $T^{\ac}_\Gamma$ satisfies the Mordell-Lang condition. Let $(\calK,\calG) \models T^{\ac}_\Gamma$ and suppose that $(\calK,\calG)$ is $\kappa$-saturated for $\kappa > |T^2|$. 
Fix $C \subseteq K$ with $|C|< \kappa$ and fix a non-algebraic unary $\calLB(C)$-type $q(x)$ and a unary $\calLA(G\cap C)$-type $p(x)$ such that $q \models \qf(p)$. By Lemma \ref{Mannalg}, $p$ must be nonalgebraic. Let $\psi(x)$ be an $\calLA(G\cap C)$-formula in $p(x)$ and let $\phi(x)$ be an $\calLB(C)$-formula in $q(x)$. By saturation, we need only show that there is an element in $G$ satisfying $\psi_A$ and $\phi$, but this follows since $\psi_A(\calG)$ is infinite and $\phi(\calK)$ is cofinite (since $\TB$ is strongly minimal). For the codensity condition, observe that by \cite[Lemma 2.2(2)]{DG06} the set $K\setminus \acl (C \cup A)$ is infinite. Thus, the codensity condition also follows from saturation of $(\cal K,\calG)$ and strong minimality of $\TB$.
\end{proof}
\end{prop}

\subsection{$\calA$-small sets}
Let $(\calB,\calA)\models T^2$. In this subsection we study $\calA$-small sets.

\begin{defn} A set $X \subseteq B$ is \textbf{$\calA$-small} if there is an $\calLB(B)$-formula $\phi(\vec{x},y)$ such that $\calB \models \forall \vec{x}\ \exists ^{<\infty}y \phi (  \vec{x}, y)$, and
$$X \subseteq \big\{b \in B: \calB \models \phi( \vec{a}, b) \text{ for some }\vec{a} \in A^{|\vec{x}|} \big\}.$$
\end{defn}

\noindent Note that a finite union of $\calA$-small sets is $\calA$-small. Moreover, $\TB$ has definable Skolem functions, then $X$ is $\calA$-small if and only if $X \subseteq f(A^n)$ for some $\calLB(B)$-definable function $f:B^n \to B$ (this follows from an easy coding argument). If $(\calB,\calA)$ is a dense pair, then the $\calA$-small sets are exactly the $\calA$-small sets in the sense of \cite{vdD98}. If $X$ is not $\calA$-small, then even if $X\subseteq \acl(A)$, this is not witnessed by finitely many formulas.

\begin{lem}
\label{notsmall}
If the pair $(\calB,\calA)$ is $\kappa$-saturated where $\kappa>|T^2|$ and if $B$ is not $\calA$-small, then $B \not \subseteq \acl ( A \cup C)$ for any $C \subseteq B$ with $|C| < \kappa$. In particular, any basis for $\calB$ over $A$ (with respect to the pregeometry induced by $\acl$) must have cardinality at least $\kappa$.
\begin{proof}
Let $C \subseteq B$ with $|C|<\kappa$ and let $\Gamma(y)$ be the partial type consisting of formulas of the form $\forall \vec{x} \left(A(\vec{x}) \rightarrow \neg\phi(\vec{x},y)\right)$ where $\phi(\vec{x},y)$ is an $(n+1)$-ary $\calL_\beta(C)$-formula such that $\phi(\vec{a},y)$ is algebraic for all $\vec{a} \in A^n$. By assumption, $\Gamma(y)$ is realizable, hence realized by some element $b \in \calB$. This $b$ is then algebraically independent over $A\cup C$.
\end{proof}
\end{lem}

\noindent For any theory $T \supseteq T^2$ which satisfies the codensity condition and for any $(\calB,\calA) \models T$, it is immediate that the set $B$ is not $\calA$-small. We have a partial converse:
\begin{lem}
\label{expODAG}
Suppose that $\TB$ is an o-minimal theory extending the theory of ordered divisible abelian groups. Let $T \supseteq T^2$ be a theory such that in every model $(\calB,\calA)\models T$ the structure $\calA$ expands a dense subgroup of $\calB$ and $B$ is not $\calA$-small. Then $T$ satisfies the codensity condition.
\begin{proof}
Let $(\calB,\calA) \models T$ be $\kappa$-saturated for $\kappa > |T^2|$ and let $C\subseteq B$ with $|C|<\kappa$. We will show that $\acl(A \cup C)$ is (topologically) codense in $B$, whence the codensity condition follows by o-minimality. Let $I$ be an interval in $B$. By Lemma \ref{notsmall}, there is an element $b \in B\setminus \acl(A \cup C)$. By density of $A$ in $B$ there is $a \in A \cap (I+b)$. But then $a- b \in I \cap B \setminus \acl(A \cup C)$.
\end{proof}
\end{lem}

\section{The Back-and-Forth System}\label{s:bnf}

Throughout this section, let $T \supseteq T^2$ be a consistent ML-theory. Let $( \calB^*_1, \calA^*_1)$ and $( \calB^*_2, \calA^*_2)$ be two $\kappa$-saturated models of $T$, where $\kappa > |T^2 | $.

\begin{asmp}\label{defn:I} Let $\calI$ be the set of all partial $\calLB$-elementary maps $\iota: B_1 \to B_2$ between finite subsets $B_1 \subseteq B_1^*$ and $B_2 \subseteq B_2^*$ such that
\begin{enumerate}
\item $B_1 \forkindep_{A_1} A_1^*$ and $B_2 \forkindep_{A_2} A_2^*$,
\item $\iota(A_1)=A_2$,
\item the restriction $\iota$ to $A_1$ is a partial $\calLA$-elementary map between $A_1$ and $A_2$,
\end{enumerate}
where $A_1 := A\cap B_1$ and $A_2 := A\cap B_2$.
\end{asmp}

\noindent Note that each map $\iota \in \calI$ is a partial $\calL^2$-isomorphism. One easily verifies the following lemma:
\begin{lem}
\label{easyremark}
Let $\iota : B_1 \to B_2$ be in $\calI$ and let $a_1\in A_1^*$ and $a_2\in A_2^*$ be such that $\iota\tp_{\calLB}(a_1|B_1) = \tp_{\calLB}(a_2|B_2)$ and $\iota\tp_{\calLA}(a_1|A_1) = \tp_{\calLA}(a_2|A_2)$. Then $\hat{\iota} := \iota\cup \big\{(a_1,a_2)\big\}$ is in $\calI$.
\end{lem}

\noindent For our pivotal result, we show that the collection $\calI$ is a ``back-and-forth system.''
A back-and-forth system (in some language) is a collection of partial isomorphisms $f:\calM\rightharpoonup \calN$ where $\calM$  and $\calN$ are ``special''  (e.g. countable, saturated)  models of some theory such that the following holds:
\begin{enumerate}
\item If $f$ is in the back-and-forth system and $a \in M$, then there is a map $g$ in the back-and-forth system which extends $f$ and which includes $a$ in its domain.
\item If $f$ is in the back-and-forth system and $b \in N$, then there is a map $g$ in the back-and-forth system which extends $f$ and which includes $b$ in its range.
\end{enumerate}
Any map in a back-and-forth system is an elementary map. The consequences of the existence of a back-and-forth system depend on the assumptions made on the domain and range of the maps in the system, but the back-and-forth method is commonly used to show completeness, countable categoricity, or some sort of quantifier reduction. 
The key fact about back-and-forth systems is that the domain and range of any map in the system have the same type.
Our back-and-forth system is modeled after the systems in \cite{vdD98} and \cite{DG06} among others, and we use it to get similar quantifier reduction results.

\begin{thm}
\label{bnf}
The set $\calI$ is a back-and-forth system in the language $\calL^2$. Therefore, each map $\iota\in \calI$ is $\calL^2$-elementary.
\begin{proof}
Let $\iota : B_1 \to B_2 \in \calI$ and $b_1 \in  B_1^*$. By symmetry it is enough to show that if $b_1\notin B_1$, then we can find $\iota'\in \calI$ extending $\iota$ such that $b_1$ is in the domain of $\iota'$. From now on, assume that $b_1\notin B_1$.\newline

\noindent\textbf{Case I. $\boldsymbol{b_1 \in B_1^* \setminus \acl(B_1 \cup A_1^*)}$:} Let $q$ be the $\calLB(B_1)$-type of $b_1$. As $\iota$ is a partial $\calLB$-elementary map, $\iota q$ is realizable in $\calB_2^*$. By the codensity condition, Definition \ref{defn:ML} (3b), we can find a realization $\iota q$ that is not in $\acl(B_2 \cup A_2^*)$. We extend $\iota$ to $\iota':B_1 \cup \{b_1\} \to B_2 \cup \{b_2\}$ by mapping $b_1$ to $b_2$. By construction, $\iota'$ is a partial $\calLB$-elementary map. It follows easily from the $\acl$-independence of $b_1$ over $B_1 \cup A_1^* $ that $\iota'$ satisfies conditions (1)-(3) of Definition \ref{defn:I}. Thus $\iota' \in \calI$.\newline

\noindent\textbf{Case II. $\boldsymbol{b_1 \in A_1^*}$:} By Lemma \ref{easyremark}, it suffices to find an element $b_2 \in A_2^*$ with $\iota\tp_{\calLB}(b_1|B_1) = \tp_{\calLB}(b_2|B_2)$ and $\iota\tp_{\calLA}(b_1|A_1) = \tp_{\calLA}(b_2|A_2)$. We consider two subcases:
\begin{description}
\item[(a) $\boldsymbol{b_1 \in \acl(B_1)}$] Let $\vec{b}$ be a tuple enumerating $B_1$ (so $\vec{b}_\alpha$ enumerates $A_1$). Set $q := \tp_{\calLB}(\vec{b})$, set $p:= \tp_{\calLA}(\vec{b}_\alpha)$, let $\phi(\vec{b},y)$ isolate the type $\tp_{\calLB}(b_1|\vec{b})$ and let $\psi(\vec{b}_\alpha,y)$ be an arbitrary formula in $\tp_{\calLA}(b_1|\vec{b}_\alpha)$. Then $(p,q,\phi,\psi)$ is a Mordell-Lang challenge for $T$ and $\big((\calB_1^*,\calA_1^*),\vec{b},b_1\big)$ is a solution. By the $\calLB$- and $\calLA$-elementarity of $\iota$ (Assumption \ref{defn:I})  the tuple $\big((\calB_2^*,\calA_2^*),\iota\vec{b}\big)$ is a contender, so by the Mordell-Lang condition, Definition \ref{defn:AE}, there is $a \in A_2^*$ such that $\big((\calB_2^*,\calA_2^*),\iota\vec{b},a\big)$ is also a solution. As $\psi$ is arbitrary, we use saturation to find some element $b_2 \in B_2^*$ which realizes $\iota\tp_{\calLA}(b_1|\vec{b}_\alpha)_A$ and such that $\calB_2^* \models \phi(\iota\vec{b},b_2)$. Now use that $\iota\vec{b}$ is an enumeration of $B_2$ and that $\phi(\iota\vec{b},y)$ isolates $\iota\tp_{\calLB}(b_1|\vec{b})$.

\item[(b) $\boldsymbol{b_1 \not \in \acl(B_1)}$] 
Let $q$ be the $\calLB(B_1)$-type of $b_1$ and let $p$ be the $\calLA(A_1)$-type of $b_1$. Then $\qf(p|_{\calL})$ is just $\qftp_{\calL}(b_1|A_1)$, so $q \models \qf (p|_{\calL})$ and $\iota q \models \iota \qf ( p|_{\calL})$. Again by the $\calLB$- and $\calLA$-elementarity of $\iota$, the type $\iota q$ is non-algebraic and $\iota \qf ( p|_{\calL}) = \qf(\iota p|_{\calL})$, so by applying the density condition, Definition \ref{defn:ML} (3a), to $( \calB^*_2, \calA^*_2)$, we find a realization $b_2 \in A_2^*$ of both $\iota p_A$ and $\iota q$.
\end{description}

\noindent\textbf{Case III. $\boldsymbol{b_1 \in \acl(B_1 \cup A_1^*) \setminus A_1^*}$:} We again consider two subcases:
\begin{description}
\item[(a) $\boldsymbol{b_1 \in \acl (B_1)}$] Let $\vec{b}$ be a tuple enumerating $B_1$, let $q(x)$ be the $\calLB(B_1)$-type of $b_1$, and let $\phi(x,\vec{y})$ be an $\calLB$-formula such that $\phi(x,\vec{b})$ isolates $q$. As $\iota$ is a partial $\calLB$-elementary map, we get $\big|\phi(\calB_1^*,\vec{b})\big|=\big|\phi(\calB_2^*,\iota\vec{b})\big|$. We claim that there is an element in $B_2^* \setminus A_2^*$ which satisfies $\phi\big(x,\iota(\vec{b})\big)$. Suppose not, so $\phi\big(\calB_2^*,\iota(\vec{b})\big) \subseteq A_2^*$. We use Case II to extend $\iota^{-1}$ to a map whose inverse is in $\calI$ including $\phi\big(\calB_2^*,\iota(\vec{b})\big)$ in its domain. This is a contradiction, as such a map would send an element in $A_2^*$ to $b_1$, which is not in $A_1^*$. Therefore, we can find an element $b_2$ in $B_2^* \setminus A_2^*$ satisfying $\phi\big(x,\iota(\vec{b})\big)$ and we extend $\iota$ by mapping $b_1$ to $b_2$. By construction, conditions (2) and (3) hold for $\iota '$ and, as $b_1 \in \acl(B_1 \cup A_1^*)$, we can easily check that (1) holds for $\iota'$ as well.

\item[(b) $\boldsymbol{b_1 \not \in \acl (B_1)}$] Take $a_1,\ldots,a_n \in A_1^*$ with $b_1 \in \acl\big(B_1\cup \{a_1,\ldots,a_n\}\big)$. By applying Case II to $a_1,\ldots,a_n$ we find a map $\iota' \in \calI$ extending $\iota$ with domain $B_1' := B_1 \cup \{a_1,\ldots,a_n\}$. Then $b_1 \in \acl(B_1')$ and so we apply the previous subcase.\qedhere
\end{description}
\end{proof}
\end{thm}

\begin{cor}\label{complete}
$T$ is complete if and only if $T_\alpha$ and $T_\beta$ are complete.
\begin{proof}
 If $T_\alpha$ and $T_\beta$ are complete, then the empty map is in the back-and-forth system constructed above. Thus, all $\kappa$-saturated models of $T$ are elementarily equivalent (as the empty map between any two such structures is $\calL^2$-elementary) so $T$ is complete.
\end{proof}
\end{cor}

\begin{cor}
\label{freeness}
Let $( \calB' , \calA')\subseteq (\calB, \calA)$ be models of $T$. Then $( \calB' , \calA' ) \preccurlyeq ( \calB , \calA )$ if and only if $\calB' \preccurlyeq \calB$, $\calA' \preccurlyeq \calA$, and $B'$ and $A$ are independent over $A'$.
\begin{proof}
$( \implies )$ Suppose $(\calB', \calA') \subseteq ( \calB , \calA )$ is an elementary substructure, and suppose that $X \subseteq B'$ is  not $\acl$-independent over $A$. Then there is an $\calLB$-formula $\phi ( \vec{x}, \vec{y} )$ such that for some $\vec{b} \in X^n$ and some $\vec{a} \in A^m$, we have $\calB \models \phi(\vec{b}, \vec{a}) \land \exists^{<\infty}x \phi (b_1, \ldots ,b_{n-1},x, \vec{a} )$.
We conclude that
$$(\calB, \calA) \models \exists \vec{y} \big(A( \vec{y} ) \land \phi(\vec{b}, \vec{y}) \land \exists^{<\infty}x \phi (b_1, \ldots ,b_{n-1},x, \vec{y} )\big).$$
By elementarity, $(\calB' , \calA')$ models this sentence as well, so $X$ is not $\acl$-independent over $A'$ either.

$( \impliedby )$ By passing to an elementary extension of $( \calB , \calA)$ if necessary, we may assume $( \calB, \calA )$ is $\kappa$-saturated.
Let $( \calB^*, \calA^*)$ be a $\kappa$-saturated elementary extension of $( \calB', \calA')$, so $B'$ and $A^*$ are independent over $A'$ by the forwards direction. Let $\calI$ be the back-and-forth system in Definition \ref{defn:I} between $( \calB , \calA)$ and $( \calB^* , \calA^*)$ and let $\vec{b}$ be a tuple in $B'$.  By Lemma \ref{finitefree}, we may assume that $\vec{b}\forkindep_{\vec{b}_\alpha}A$ and $\vec{b}\forkindep_{\vec{b}_\alpha}A^*$, so the identity map on $\vec{b}$ is a map in $\calI$. By Theorem \ref{bnf}, we conclude that  type which $\vec{b}$ realizes 
in $( \calB^*, \calA^*)$ is the same as the type it realizes 
in $( \calB , \calA )$.
Since $( \calB', \calA') \preccurlyeq ( \calB ^*, \calA ^*)$ we get that $( \calB', \calA') \preccurlyeq ( \calB, \calA)$ as well.
\end{proof}
\end{cor}

\begin{defn} An $\calL^2(B)$-formula is called \textbf{special} if it is of the form
$\theta ( \vec{y} ) = \exists \vec{x} \left(A ( \vec{x} ) \wedge \psi_A( \vec{x}) \wedge \phi ( \vec{x} , \vec{y} )\right)$ where $\phi ( \vec{x} , \vec{y} )$ is an $\calLB$-formula, and $\psi ( \vec{x} )$ is an $\calLA$-formula.
\end{defn}

\begin{thm}
\label{special}
Every $\calL^2(B)$-formula is equivalent in $T$ to a boolean combination of special formulas.
\begin{proof}
By removing and re-introducing parameters, it is enough to show that this is true for $\calL^2$-formulas without parameters.
Let $( \calB , \calA )$ be a $\kappa$-saturated model of $T$ where $\kappa > |T^2|$ and let $\calI$ be the back-and-forth system in Definition \ref{defn:I} between $( \calB , \calA)$ and itself.
Let $\vec{b} = (b_1, \ldots ,b_n)$ and $\vec{d} = (d_1, \ldots ,d_n)$ be tuples from $B$ that satisfy the same special formulas. It suffices to show that $\vec{b}$ realizes the same $\calL^2$-type as $\vec{d}$. For this, it is enough to find $\iota \in \calI$ that sends $\vec{b}$ to $\vec{d}$.

Let $r \leq n$ be the $\acl$-rank of $(b_1, \ldots ,b_n)$ over $A$.
Without loss of generality, we may assume $b_1, \ldots , b_r$ are $\acl$-independent over $A$.
Then for every $\calLB$-formula $\phi(\vec{x},\vec{y})$ and each $i \in \{1,\ldots,r\}$, we must have
$$(\calB,\calA) \models \neg \exists \vec{x} \left(A(\vec{x}) \land \phi (\vec{x},b_1, \ldots , b_r) \land \exists ^{< \infty} z \phi (\vec{x},b_1, \ldots ,b_{i-1}, z,b_{i+1},\ldots,b_r)\right).$$
By the assumption that $\vec{b}$ and $\vec{d}$ satisfy the same special formulas we conclude that $d_1, \ldots ,d_r$ are also $\acl$-independent over $A$.

For each $i \in \{ r+1, \ldots ,n \}$, set $B_i := A \cup \{b_1,\ldots,b_{i-1}\}$ and let $\phi_i(\vec{a},b_1,\ldots,b_{i-1},z)$ be a $\calLB(B_i)$-formula isolating the type of $b_i$ over $B_i$ where $\vec{a}$ is a tuple in $A$ (by adding dummy variables, we may assume that $\vec{a}$ is the same for each formula). We want to find a tuple $\vec{c} \in A^{|\vec{a}|}$ such that $\tp_{\calLA}( \vec{a}) = \tp _{\calLA} ( \vec{c})$, $\tp_{\calLB}( \vec{a} ) = \tp _{\calLB} ( \vec{c})$, and $\calB \models \phi _i (\vec{c} , d_1, \ldots ,d_{i-1}, d_i)$ for each $i \in \{ r+1, \ldots ,n \}$. Fix $\psi(\vec{x}) \in \tp_{\calLA}(\vec{a})$ and $\phi(\vec{x}) \in \tp_{\calLB}(\vec{a})$. Note that the formula
$$\theta(\vec{y}) :=\exists \vec{x} \big(A(\vec{x}) \wedge \psi_A(\vec{x}) \wedge \phi(\vec{x}) \wedge \bigwedge_{i=r+1}^n \phi_i(\vec{x}, y_1,\ldots,y_i)\big)$$
is a special formula and that $(\calB,\calA) \models \theta(\vec{b})$. Thererefore, $(\calB,\calA) \models \theta(\vec{d})$ and so, by saturation, we find a tuple $\vec{c}$ with the desired properties.

By repeated application of Lemma \ref{easyremark} and the fact that the empty map is in $\calI$, there is a map $\iota \in \calI$ sending $\vec{a}$ to $\vec{c}$. Proceeding as in Case I of Theorem \ref{bnf}, we extend $ \iota$ to a map $\iota' \in \calI$ sending also $\{b_1,\ldots,b_r\}$ to $\{d_1,\ldots,d_r\}$. Finally, we extend $\iota'$ to a map $\iota'' \in \calI$ sending $\{b_{r+1},\ldots,b_n\}$ to $\{d_{r+1},\ldots,d_n\}$ recursively: if $b_i \in A$ for $i = r+1,\ldots,n$, then $b_i$ must be a component of $\vec{a}$ since $\phi_i$ isolates the type of $b_i$ over $B_i \supseteq A$. If $b_i \not\in A$ then by the argument in Case III of Theorem \ref{bnf} and since $\phi_i$ isolates the type of $b_i$, we can extend by sending $b_i$ to $d_i$.
\end{proof}
\end{thm}

\noindent A theory is said to be \textbf{near model complete} if every formula is equivalent to a boolean combination of existential formulas. The following is immediate from Theorem \ref{special}:
\begin{cor} \label{nmc}
If $T_\beta$ and $T_\alpha$ are model-complete, then $T$ is near model complete.
\end{cor}

\noindent As a remark, a theory can be near model complete but not model complete. A proof is given in \cite{Ro59} that the theory of the pair $(\overline{\R}, \calA)$, where $\overline{\R}$ is the real field and $\calA$ is the field of real algebraic numbers, is not model complete (in the language of ordered rings with an additional unary predicate).

\section{Types, Open core, and NIP}\label{s:opencore}
Let $T \supseteq T^2$ be a consistent ML-theory and let $(\calB,\calA)$ be a $\kappa$-saturated model of $T$ where $\kappa > |T^2|$. In this section, we prove two important preservation results. The first result states that if $\TB$ is equipped with a definable topology satisfying certain weak conditions, then every open subset of $B^n$ definable in $(\calB,\calA)$ is already definable in $\calB$. Thus expanding $\calB$ by $\calA$ does not introduce any new open sets. The second result concerns the preservation of model-theoretic tameness: if $\TB$ and $\TA$ are both complete NIP theories, then $T$ is NIP as well. Before we can prove these theorems, we have to study types in ML-pairs in more detail.

\subsection{Types}
In this subsection, we use the back-and-forth system constructed in the previous section to characterize some $\calL^2$-types. For the remainder of this subsection, let $C$ be a finite subset of $B$ such that $C \forkindep_{A\cap C}A$ and let $\calI$ be the back-and-forth system in Definition \ref{defn:I} between $(\calB,\calA)$ and itself.
\begin{lem}
\label{atype} Let $\vec{a}_1,\vec{a}_2 \in A^n$ be such that
\begin{enumerate}
\item $\tp_{\calLB}(\vec{a}_1|C) = \tp_{\calLB}(\vec{a}_2|C)$, and
\item $\tp_{\calLA}(\vec{a}_1|A\cap C) = \tp_{\calLA}(\vec{a}_2|A\cap C)$.
\end{enumerate}
Then $\tp_{\calL^2}(\vec{a}_1|C) = \tp_{\calL^2}(\vec{a}_2|C)$.
\begin{proof}
The identity map on $C$ is in $\calI$ and by repeated application of Lemma \ref{easyremark}, the map $\iota:C\vec{a}_1 \to C\vec{a}_2$ which is the identity on $C$ and sends $\vec{a}_1$ to $\vec{a}_2$ is in $\calI$. Thus $\vec{a}_1$ and $\vec{a}_2$ have the same $\calL^2(C)$-type by Theorem \ref{bnf}.
\end{proof}
\end{lem}

\noindent From the above lemma we conclude that $\calL^2$-definable subsets of $A$ are determined by $\calLB$-definable and $\calLA$-definable sets in the following way:

\begin{cor}
\label{acor}
Every $\calL^2(C)$-definable subset $X \subseteq A^n$ is a boolean combination of $\calLB(C)$-definable subsets of $B^n$ and $\calLA(A\cap C)$-definable subsets of $A^n$.
\end{cor}

\begin{defn}\label{Dindep}
 For $n\in \N$, we define $D_n(C)$ to be the set
\[
\{\vec{x} \in B^n: \vec{x} \text{ is } \acl \text{-independent over }A\cup C\}.
\]
\end{defn}

\begin{lem}
\label{Dlemma}
Let $\vec{d}_1,\vec{d}_2 \in D_n(C)$ be such that $\tp_{\calLB}(\vec{d}_1|C) = \tp_{\calLB}(\vec{d}_2|C)$. Then $\tp_{\calL^2}(\vec{d}_1|C) = \tp_{\calL^2}(\vec{d}_2|C)$.
\begin{proof}
Again it suffices to show the statement of the lemma for every finite subset of $C$. Therefore, by Lemma \ref{finitefree}, we may assume that $C$ is finite. The identity map on $C$ is in $\calI$, so let $\iota:C\vec{d}_1 \to C\vec{d}_2$ be the extension of the identity map on $C$ and sends $\vec{d}_1$ to $\vec{d}_2$. We will now show that $\iota \in \calI$. By assumption $\iota$ is a partial $\calLB$-elementary map. Since $\vec{d}_1,\vec{d}_2\in D_n(C)$, we easily get that
\begin{enumerate}
\item $A\cap (C\vec{d_1})=A\cap C=A\cap (C\vec{d_2})$,
\item $C\vec{d}_1\forkindep_{A\cap C}A$, $C\vec{d}_2 \forkindep_{A\cap C}A$.
\end{enumerate}
Thus the restriction of $\iota$ to $A\cap(C\vec{d}_1)$ is $\calLA$-elementary, so $\iota\in \calI$ and $\vec{d}_1$ and $\vec{d}_2$ have the same $\calL^2(C)$-type.
\end{proof}
\end{lem}

\subsection{Open sets}
In this subsection, we suppose that $\TB$ is equipped with a \textbf{definable topology}, that is, there is $n>0$ and a distinguished $(1+n)$-ary $\calLB$-formula $\tau(x,\vec{y})$ such that for every model $\calB\models \TB$, the family of definable sets
$$\big\{\tau(\calB,\vec{d}):\vec{d} \in B^n\big\}$$
forms a basis for a topology on $\calB$. For each $m$ and each $\vec{d} = (\vec{d}_1,\ldots,\vec{d}_m) \in B^{n\times m}$, we let
$$U_{\vec{d}}:=\big\{(x_1,\ldots,x_m) \in B^m: x_i \in \tau(\calB,\vec{d_i})\text{ for each } i = 1,\ldots,m\big\}.$$

\noindent We assume that $\tau(\calB,\vec{d})$ is either empty or infinite for every $\vec{d} \in B^n$. We also assume that for every open set $V \subseteq B^m$ and every $\vec{x} \in V$, the set
$$\{\vec{d} \in B^{n\times m}:\vec{x} \in U_{\vec{d}}\text{ and }U_{\vec{d}}\subseteq V\}$$
has non-empty interior in $B^{n\times m}$. This second assumption is Assumption (I) in Boxall and Hieronymi \cite{BH12}. These assumptions are satisfied when $\calB$ is an o-minimal structure or a $p$-adically closed field. The main theorem for this section describes the open sets definable in $(\calB,\calA)$.

\begin{thm}
\label{opencore}
Every open set definable with parameters in $(\calB,\calA)$ is definable (perhaps with additional parameters) in $\calB$.
\begin{proof}
Let $X\subseteq B^m$ be open and $\calL^2$-definable with parameters from a finite set $C$. By increasing $C$, we may assume that $C\forkindep_{A\cap C}A$ (see Lemma \ref{finitefree}). We will now prove that $X$ is $\calLB(C)$-definable. By Corollary 3.1 in \cite{BH12}, it suffices to show that the set $D_m(C)$ (see Definition \ref{Dindep}) has the following properties:
\begin{enumerate}[(i)]
\item $D_m(C)$ is dense in $B^m$;
\item for every $\vec{b} \in D_m(C)$ and every open set $V\subseteq B^m$, if $\tp_{\calLB}(\vec{b}|C)$ is realized in $V$, then $\tp_{\calLB}(\vec{b}|C)$ is realized in $V \cap D_m(C)$;
\item for every $\vec{b} \in D_m(C)$, $\tp_{\calL^2}(\vec{b}|C)$ is implied by $\tp_{\calLB}(\vec{b}|C)$ and membership in $D_m(C)$.
\end{enumerate}

For property (i), fix $\vec{d}_1,\ldots,\vec{d}_m \in B^n$ such that $\tau(\calB,\vec{d}_i)$ is nonempty for each $i$. By repeatedly invoking the codensity condition, we realize a tuple $\vec{b}$ in $U_{\vec{d}}\cap D_m(C)$.

For property (ii), let $\vec{b}$ and $V$ be given and fix a realization $\vec{b}'$ of $\tp_{\calLB}(\vec{b}|C)$ in $V$. By our assumptions, the set
$$\{\vec{d} \in B^{n\times m}:\vec{b}' \in U_{\vec{d}}\text{ and }U_{\vec{d}}\subseteq V\}$$
has nonempty interior. Since nonempty open sets are assumed to be infinite, we can find $\vec{d}$ in this set such that $\vec{d}$ is $\acl$-independent over $C\vec{b}'$. Since $\vec{b}'$ is $\acl$-independent over $C$ and since $\acl$ is a pregeometry, $\vec{b}'$ must be $\acl$-independent over $C\vec{d}$. By repeatedly invoking the codensity condition, we find a tuple $\vec{b}''$ realizing $\tp_{\calLB}(\vec{b}'|C\vec{d})$ with $\vec{b}'' \in D_m(C)$. In particular, $\vec{b}''$ is in $U_{\vec{d}}$.

Property (iii) is just Lemma \ref{Dlemma}.
\end{proof}
\end{thm}

\subsection{NIP for ML-theories}
In this subsection, we show that if both $\TB$ and $\TA$ are complete NIP theories then $T$ is NIP. To do this, we apply a result of Chernikov and Simon. We restate a version of this result as Fact \ref{CS} so that it applies more directly to our case.

\begin{defn}
Let $\tilde{T} \supseteq T^2$ be a complete $\calL^2$-theory, let $\theta(\vec{x},\vec{y})$ be an $\calL^2$-formula and let $(\calB,\calA)$ be a $\kappa$-saturated model of $\tilde{T}$ for $\kappa >|T^2|$.
\begin{enumerate}
\item $\theta$ is said to be \textbf{NIP} if there is no $\calL^2$-indiscernible sequence $(\vec{a}_i)_{i \in \omega}$ from $B^{|\vec{x}|}$ and no $\vec{b} \in B^{|\vec{y}|}$ such that $(\calB,\calA) \models \theta(\vec{a}_i,\vec{b})$ if and only if $i$ is odd.
\item $\tilde{T}$ is said to be \textbf{NIP} if every $\calL^2$-formula is NIP.
\item $\theta$ is said to be \textbf{NIP over $A$} if there is no $\calL^2$-indiscernible sequence $(\vec{a}_i)_{i \in \omega}$ from $A^{|\vec{x}|}$ and no $\vec{b} \in B^{|\vec{y}|}$ such that $(\calB,\calA) \models \theta(\vec{a}_i,\vec{b})$ if and only if $i$ is odd.
\item $\tilde{T}$ is said to be \textbf{NIP over $A$} if every $\calL^2$-formula is NIP over $A$.
\end{enumerate}
\end{defn}

\begin{fact}\cite[Theorem 2.4]{CS13}
\label{CS}
Let $\tilde{T} \supseteq T^2$ be a complete $\calL^2$-theory, let $\theta(\vec{x}\vec{y},\vec{z})$ be an $\calL^2$-formula and let $(\calB,\calA)$ be a $\kappa$-saturated model of $\tilde{T}$ for $\kappa >|T^2|$. If $\theta$ is NIP and if $\tilde{T}$ is NIP over $A$ then $\exists \vec{x}\big(A(\vec{x}) \wedge\theta(\vec{x}\vec{y},\vec{z})\big)$ is NIP.
\end{fact}

\begin{thm}\label{NIP}
If both $\TB$ and $\TA$ are complete NIP theories then so is $T$.
\begin{proof}
As NIP formulas are preserved by boolean operations, and as $\TB$ and $\TA$ are NIP, we see from Corollary \ref{acor} that $T$ is NIP over $A$. By Theorem \ref{special}, it suffices to show that each $\calL^2$-formula of the form
$$\theta(\vec{y},\vec{z}) = \exists \vec{x} \big(A ( \vec{x} ) \wedge \psi_A( \vec{x}) \wedge \phi ( \vec{x} \vec{y},\vec{z} )\big)$$
is NIP, where $\phi ( \vec{x}\vec{y} , \vec{z} )$ is an $\calLB$-formula, and $\psi ( \vec{x} )$ is an $\calLA$-formula. However, this follows from Fact \ref{CS}, noting that $ \psi_A( \vec{x}) \wedge \phi ( \vec{x} \vec{y},\vec{z} )$ is NIP.
\end{proof}
\end{thm}

\section{Pairs of distinct o-minimal structures and ordered vector spaces}\label{s:oml}

Let $\TB$ and $\TA$ be o-minimal theories extending the theory of dense linear orders without endpoints and suppose that $\calL\supseteq \{<\}$. In this special case, we call an ML-theory $T \supseteq T^2$ an \textbf{o-ML-theory} and we call a model $(\calB,\calA) \models T$ an \textbf{o-ML-pair}. One particular example of an o-ML-pairs is a dense pair of o-minimal structures as studied in \cite{vdD98}, which serves as the inspiration for the definition of the broader class of o-ML-pairs.

\subsection{Properties of o-ML-theories}
As is proven in Lemma \ref{expODAG}, if $\TB$ and $\TA$ extend the theory of ordered abelian groups, $(\calB,\calA) \models T^2$, $\calA$ is a dense subgroup of $\calB$, and $\calB$ is not $\calA$-small, then the theory of the pair $(\calB,\calA)$ satisfies the codensity condition.  
Now that we are in the o-minimal setting, we make liberal use of the order topology.  We exploit the fact that the density condition is related to topological density in the following way:

\begin{lem}\label{omldensity}
Suppose that $\calLA \subseteq \calLB$ and that $\TA$ admits quantifier elimination. Let $T \supseteq T^2$ and suppose that for every model $(\calB,\calA) \models T$, the topological closure of $A$ in $B$ is $\calLA$-definable without parameters. Then $T$ satisfies the density condition.
\begin{proof}
Fix $\kappa > |T^2|$, a $\kappa$-saturated model $(\calB,\calA) \models T$, a subset $C \subseteq B$ with $|C|<\kappa$, and an $\calLB(C)$-type $q$. 
Let $p$ be any $\calLA(A\cap C)$-type such that $q \models \qf (p|_{\calL})$, so $q \models \qf(p)$. By quantifier elimination for $\TA$, the type $p$ is completely determined by its quantifier-free part, so for $a \in A$, if $(\calB,\calA) \models \qf(p)(a)$, then $(\calB,\calA) \models \qf(p_A)(a)$ and so $(\calB,\calA) \models p_A(a)$. Therefore, it suffices to find $a \in A$ realizing $q$. 
By assumption, the closure of $A$ in $B$ is a finite union of $\calLA$-definable points and open intervals. 
Since $p$ is non-algebraic, one (and hence all) realizations of $p$ are in one of these open intervals; call it $I$. 
For each $\phi(x) \in q(x)$, we may assume that $\phi(x)$ defines an interval contained in $I$ and so by density of $A$ in $I$, we have that $(\calB,\calA) \models \exists x \big(A(x) \wedge \phi(x)\big)$. We are done by saturation.
\end{proof}
\end{lem}

\noindent Though the conditions in the lemma above may seem somewhat peculiar, the fact that we do not assume that $A$ is dense in $B$ gives us some additional flexibility, as we will see in the following example.

\subsubsection*{Real closed fields with a Mann subgroup}
Let $\Gamma$ be a dense, divisible, multiplicative subgroup of $\R^{>0}$ with the Mann property (see Definition \ref{def:mann}). We axiomatize the pair $(\R,\Gamma)$ as follows: set $\calLA:= \big\{0,1, \cdot,<, (\gamma)_{\gamma \in \Gamma} \big\}$ and let $\TA$ be the $\calLA$-theory of $\Gamma\cup \{0\}$ (so $\TA$ is the theory of ordered divisible abelian groups with distinguished elements and a point at $-\infty$). Set $\calLB :=\big\{0,1, \cdot,+,- ,<,(\gamma)_{\gamma \in \Gamma}\big\}$ and let $\TB$ be the $\calLB$-theory of $\R$. We let $T^{\rc}_\Gamma \supseteq T^2$ be the theory stating that for $(\calR,\calG)\models T^{\rc}_\Gamma$ and for every $\Q$-linear equation $\sum_{i = 1}^n q_ix_i = 1$, each non-degenerate solution in $\calG$ is among one of the solutions in $\Gamma$.

\begin{prop}\label{Mannoml}
$T^{\rc}_\Gamma$ is an o-ML-pair.
\begin{proof}
We first show that $T^{\rc}_\Gamma$ satisfies the Mordell-Lang condition. Let $(p,q,\varphi,\psi)$ be a solvable Mordell-Lang challenge for $T^{\rc}_\Gamma$ with solution $\big((\calR, \calG), \vec{c}, a \big)$. 
Since we assume $a\in G$ is algebraic over $\vec{c}$ and that $\vec{c}\forkindep_{\vec{c}_\alpha}G$, we have that $\vec{c}_\alpha$ and $a$ are algebraically dependent over $\F(\Gamma)$. By Lemma 5.12 in \cite{DG06}, we have that $\vec{c}_\alpha$ and $a$ are multiplicatively dependent over $\Gamma$, so we have 
\[
a = c_1^{p_1}\ldots c_n^{p_n}\gamma_1^{q_1}\ldots\gamma_m^{q_m}
\]
where $(c_1,\ldots,c_n)$ is an enumeration of $\vec{c}_\alpha$, where $\gamma_1,\ldots,\gamma_m \in \Gamma$, and where each $p_i$ and each $q_j$ is a rational number. Let $\theta(\vec{x}_\alpha,y)$ be the formula $x_1^{p_1}\ldots x_n^{p_n}\gamma_1^{q_1}\ldots\gamma_m^{q_m}=y$. 
Then $(\calR, \calG) \models \theta(\vec{c}_{\alpha},a)$, and since $\theta(\vec{x}_\alpha,y)$ isolates the $\calLB$- and $\calLA$-type of $a$ over $\vec{c}$, we conclude that
\[
p(\vec{x}) \vdash \forall y\big( \theta(\vec{x}_\alpha,y) \rightarrow \phi(\vec{x},y)\big) \text{ and } q(\vec{x}_\alpha) \vdash \forall y\big( \theta(\vec{x}_\alpha,y) \rightarrow \psi(\vec{x}_\alpha,y)\big).
\]
Thus, any contender $\big((\calR',\calG'),\vec{d}\ \big)$ to the Mordell-Lang challenge can be extended to a solution $\big((\calR',\calG'),\vec{d},a' \big)$ by setting $a':= d_1^{p_1}\ldots d_n^{p_n}\gamma_1^{q_1}\ldots\gamma_m^{q_m}$ where $(d_1,\ldots,d_n)$ enumerates $\vec{d}_\alpha$.

Now let $(\calR,\calG) \models T^{\rc}_\Gamma$ and suppose that $(\calR,\calG)$ is $\kappa$-saturated for $\kappa > |T^2|$. The density condition follows from density of $G$ in $\calR^{>0}$, quantifier elimination for ordered divisible abelian groups, and Lemma \ref{omldensity}.
We deduce codensity by showing that $\acl(G \cup C)$ is codense in $R$ for any $C\subseteq R$ with $|C|<\kappa$ (the codensity condition follows by o-minimality).
Let $I\subseteq R^{>0}$ be an interval.
By Lemma 6.1 in \cite{DG06}, $R$ is not $\calG$-small, so by Lemma \ref{notsmall}, there is an element $r \in R^{>0}\setminus \acl(G \cup C)$. By density of $G$ in $R^{>0}$ there is $g \in G \cap (I\cdot r)$. But then $\frac{g}{r} \in I \cap R \setminus \acl(G \cup C)$.
\end{proof}
\end{prop}

\noindent For the remainder of this subsection, fix an o-ML-theory $T$. We list here the consequences of Theorems \ref{bnf} and \ref{special}.

\begin{cor}\label{nmcoML}
If $\TB$ and $\TA$ are both complete, then $T$ is complete as well. If in addition $\TB$ and $\TA$ are model complete, then $T$ is near model complete.
\end{cor}

\noindent We also have a characterization of the open core of an o-ML-pair by fact that the open core of $\calB$ is interdefinable with $\calB$ by o-minimality.

\begin{cor}\label{ocoML}
For an o-ML-pair $(\calB , \calA)$, the open core of $( \calB, \calA)$ is interdefinable with $\calB$ (so $\TB$ is an open core of $T$).
\end{cor}

\noindent Finally, we can conclude the following from Theorem \ref{NIP} and the fact that o-minimal theories are NIP:

\begin{cor}
If $\TA$ and $\TB$ are complete, then $T$ is NIP.
\end{cor}

\subsection{Pairs of ordered vector spaces}
In this subsection, we fix a subfield $K \subseteq \R$ with $\Q \subsetneq K$ and examine the pair $( \tilde{ \R }, \tilde{ \Q})$ where $\tilde{ \R } := \big( \R, 0, 1, <, +, (\lambda_k)_{k \in K}\big)$ is the reals as an ordered vector space over $K$, and $\tilde{ \Q } := \big( \Q, 0,1, <, +, (\lambda_q)_{q \in \Q}\big)$ is $\Q$ as an ordered vector space over itself (where $\lambda_k$ denotes the function $x \mapsto kx$). We will see that the first order theory of this pair is an o-ML-theory.\newline

\noindent Let $\calLB := \big\{ <,+,0,1, (\lambda_k)_{k \in K}\big\}$ be the language of ordered $K$-vector spaces with distiguished positive element 1 and let $\calLA \subseteq \calLB$ be the sublanguage of ordered $\Q$-vector spaces. Let $\frI$ denote the collection of all finite $\Q$-linearly independent subsets of $K$.

\begin{defn}
Let $T^d_K$ be the $\calL^2$-theory whose models $(\calR, \calQ)$ satisfy the following statements:
\begin{enumerate}
\item{$\calR$ is an ordered $K$-vector space with distinguished positive element 1.}
\item{$\calQ$ is an ordered $\Q$-vector subspace of $\calR$ with distinguished positive element 1.}
\item{$Q$ is dense in $R$.}
\item{For all $n \in \N$ and all $\{k_1,\ldots,k_n \} \in \frI$ there is $r \in R$ such that $r \not \in  \lambda_{k_1}( Q) + \ldots + \lambda_{k_n}( Q)$.}
\item{ For all $n \in \N$ and all $\{k_1, \ldots ,k_n \} \in \frI$, and for all $x_1,\dots,x_n\in Q$
\[
 \lambda_{k_1}(x_1) + \cdots +\lambda_{k_n}(x_n)=0 \Longrightarrow \bigwedge_{i=1}^n x_i= 0.\] }
\end{enumerate}
\end{defn}

\noindent Note that the structure $(\tilde{ \R}, \tilde{ \Q })$ described above is a model of this theory. Fix $(\calR,\calQ)\models T^d_K$.
The following lemma illustrates the complementary nature of how $K$ and $\calQ$ interact over $\Q$.
For the rest of this section, fix a $\Q$-linear basis $Z$ for $K$.

\begin{lem}
\label{Qindependence}
If $X \subseteq Q$ is $\Q$-linearly independent, then $X$ is $K$-linearly independent.
Moreover, for every $n$ and every $\vec{k} \in (K^{\times})^n$, there are $\vec{q}_1,\ldots,\vec{q}_m \in \Q^n$ (with $q_{i,j} \neq 0$ for some $i,j$) such that
$$T^d_K \models \forall \vec{x}  \in Q^n  \Big( \sum_{j=1}^n \lambda_{k_j} (x_j) = 0 \leftrightarrow \bigwedge_{i=1}^m\sum_{j=1}^n \lambda_{q_{i,j}} (x_j) = 0\Big).$$

\begin{proof}
We prove the ``Moreover,'' since the contrapositive of the first claim follows immediately from the second.
Take $\vec{k} \in (K^{\times})^n$ and choose $\{b_1,\ldots,b_m\} \subseteq Z$ so that we can write $k_j = \sum_{i=1}^m q_{i,j} b_i$ for $j=1, \ldots , n$ where $ q_{i,j} \in \Q$. As $k_j \neq 0$, there is an $i$ for each $j$ with $q_{i,j} \neq 0$. By linearity we obtain the following equalities:
$$\sum_{j=1}^n \lambda_{k_j} (x_j) = \sum_{j=1}^n \lambda_{\sum_{i=1}^m q_{i,j} b_i}(x_j)=\sum_{j=1}^n \sum_{i=1}^m \lambda_{ q_{i,j} b_i} (x_j) =  \sum_{i=1}^m \lambda_{ b_i}\Big( \sum_{j=1}^n \lambda_{ q_{i,j}}(x_j) \Big)$$
for all $\vec{x} \in Q^n$.
Since $b_1,\ldots,b_m$ are $\Q$-linearly independent, we know by Axiom scheme (5) that
\[
\sum_{i=1}^m \lambda_{b_i}\Big( \sum_{j=1}^n \lambda_{q_{i,j}}(x_j) \Big) = 0 \Longleftrightarrow \bigwedge_{i=1}^m \sum_{j=1}^n \lambda_{q_{i,j}}(x_j) =0.
\]
for all $\vec{x} \in Q^n$.
\end{proof}
\end{lem}

\begin{cor}
The theory $T^d_K$ is an o-ML-theory.
\begin{proof}
We first show that $T^d_K$ satisfies the Mordell-Lang condition. 
Let $(p,q,\varphi,\psi)$ be a solvable Mordell-Lang challenge for $T^d_K$ with solution $\big((\calR, \calQ), \vec{c}, a \big)$. 
Since we assume $a\in Q$ is algebraic over $\vec{c}$ and that $\vec{c}\forkindep_{\vec{c}_\alpha}Q$, it must hold that $a$ is algebraic over $\vec{c}_\alpha = (c_1, \ldots ,c_n)$. 
It follows from quantifier elimination for ordered vector spaces that any algebraic formula in the language $\calLB$ is equivalent to a positive boolean combination of linear equations of the form $\lambda_{k_0}(1)+\sum_{i=1}^m \lambda_{k_i} (x_i) = 0$, where $k_0, \ldots ,k_m \in K$. We may assume that $1$ is a component of $\vec{c}_\alpha$, so there are $k_1,\ldots,k_n \in K$ such that $\sum_{i=1}^n \lambda_{k_i} (c_i)  = a$.  
By Lemma \ref{Qindependence}, we see that there are $q_1,\ldots,q_n \in \Q$ such that $\sum_{i=1}^n\lambda_{q_i} (c_i)  = a$. 
As in the proof of Proposition \ref{Mannoml}, we see that any contender $\big((\calR',\calQ'),\vec{d}\ \big)$ to the Mordell-Lang challenge can be extended to a solution $\big((\calR',\calQ'),\vec{d},a' \big)$ by setting $a':= \sum_{i=1}^n\lambda_{q_i} (d_i)$ where $(d_1,\ldots,d_n) = \vec{d}_\alpha$.

The density condition follows from Lemma \ref{omldensity}, noting that $\calLA \subseteq \calLB$ and that $Q$ is dense in $R$.
To see that the codensity condition holds, we appeal to Lemma \ref{expODAG} and Axiom (4), which easily implies that $R$ is not $\calQ$-small in light of quantifier elimination for ordered vector spaces.
\end{proof}
\end{cor}

\noindent Since the theory of ordered vector spaces is complete, we conclude that $T^d_K$ is complete. Moreover, the theory of ordered vector spaces admits quantifier elimination, so we can deduce the following by Corollaries \ref{nmcoML} and \ref{ocoML}:

\begin{cor}\label{thma}
If $(\calR, \calQ)\models T^d_K$, then every $\calL^2(R)$-definable subset of $R^n$ is a boolean combination of sets of the form
\[
\bigcup_{\vec{q} \in Q^m} \big\{ \vec{a} \in R^n \ : \ (\vec{q},\vec{a}) \in X\big\},
\]
where $X\subseteq R^{m+n}$ is $\calLB(R)$-definable. Furthermore, every $\calL^2(R)$-definable open subset of $R^n$ is already $\calLB(R)$-definable.
\end{cor}

\noindent There is a dichotomy among the kind of models $T^d_K$ can have based on whether the dimension of $K$ over $\Q$ is finite or infinite.

\begin{cor}
If the dimension of $K$ over $\Q$ is infinite, then the structure $(K, \tilde{\Q})$ is a prime model of $T^d_K$. Moreover, $(\acl ( Q ), \calQ)$ is always an elementary substructure of $(\calR,\calQ)$.
\begin{proof}
One easily checks that $(K,\Q)$ is indeed a model of $T^d_K$ and that it canonically embeds into every other model of $T^d_K$, so it suffices to check that this embedding is elementary. By Corollary \ref{freeness} and quantifier elimination for ordered vector spaces, the substructure $(K, \Q)$ of model $( \calR, \calQ )$ is an elementary substructure if and only if $K$ and $Q$ are independent over $\Q$, and this follows since 
$\acl$-independence is the same as $K$-linear independence and 
there are no $K$-linearly independent subsets of $K$.
The ``moreover'' statement follows by Corollary \ref{freeness} as well.
\end{proof}
\end{cor}

\noindent We now characterize when $T^d_K$ is a decidable theory.

\begin{thm}
\label{decidability}
The theory $T^d_K$ is decidable if and only if $K$ has a computable presentation as an ordered field and there is a recursive algorithm for ascertaining $\Q$-linear independence for finite subsets of $K$.
\begin{proof}
With regards to the forward direction, if $K$ had no computable presentation as an ordered subfield of $\R$, then either $K$ would not be recursively enumerable or the order relation of the theory would not be decidable, hence $T^d_K$ could not be decidable. Similarly if there were no recursive algorithm for determining the $\Q$-linear independence of a given finite set of elements of $K$, it would be impossible to recursively check that a given $\calL^2$-sentence falls in the Axiom scheme (5).

For the other direction, since $T^d_K$ is a complete theory it suffices to show the axioms are recursively enumerable.
Since $K$ and $\Q$ are computable ordered fields, it follows immediately that Axioms (1) and (2) are computable. Axiom (3) is finite, hence computable.
The recursive enumerability of $\frI$ follows from the existence of a recursive algorithm for ascertaining the $\Q$-linear independence of finite subsets of $K$, so the Axiom schemes (4) and (5) are recursively enumerable.
\end{proof}
\end{thm}

\noindent We remark that a sufficient condition for having a recursive algorithm to ascertain the $\Q$-linear independence of any finite set of elements in $K$ is the existence of a computable basis for $K$ over $\Q$.
There are numerous examples of fields $K$ which are known to have a computable presentation and a computable basis as a vector space over $\Q$, including the following:

\begin{ex}
We note that it is known, as exposited in Miller \cite{Mi10}, that any field $K \supseteq \Q$ that is computably presentable and has a computable transcendence basis also has a computable $\Q$-linear basis. Thus, for the following choices of $K$, the hypotheses of Theorem \ref{decidability} are satisfied:
\begin{enumerate}
\item{
By \cite{Mi09} the field $K := \Q (\sqrt{p_1},\sqrt{p_2},\ldots)$ where $p_n$ is the $n^{th}$ prime is computably presentable, with a clear choice for computable basis. 
}

\item{The field $K := \R^{\operatorname{alg}}$ of real algebraic numbers is computably presentable.}

\item{
By \cite{Le09} the field $K := \Q(e)$ is computably presentable, with computable transcendence basis $\{ e \}$.
}

\item{By using Taylor series to expand $\pi$ it is easy to show by the methods used in \cite{Le09} that the field $K := \Q( \pi )$ is computably presentable, with computable transcendence basis $\{ \pi \}$.
}
\end{enumerate}
\end{ex}

\section{Real closed field with a predicate for a pseudo real closed subfield}\label{s:prc}
In this section, we consider a real closed field with a predicate for a dense pseudo real closed subfield with $n$ orderings where $n\geq 2$. Let
$$\calLA :=\{0,1,+,\cdot,-,<_1,<_2,\ldots,<_n\}.$$
An \textbf{$n$-ordered field} is an $\calLA$-structure $\calK = (K, \ldots)$ such that $(K,0,1,+,\cdot,-,<_i)$ is an ordered field for $i=1,\ldots,n$. Let $\TA$ be the theory of $n$-ordered fields which satisfy the following two axioms of van den Dries \cite{vdD78}:
\begin{itemize}
\item $<_i$ and $<_j$ induce different interval topologies for $1\leq i <j\leq n$;
\item for each irreducible $P(T,X) \in K[T,X]$ and $a \in K$ such that $P(a,X)$ changes sign on $K$ with respect to each ordering $<_i$, there are $c,d \in K$ with $P(c,d) = 0$.
\end{itemize}
Then $\TA$ is the model companion to the theory of $n$-ordered fields, and we say that $\calK\models \TA$ is a \textbf{pseudo real closed field}. Pseudo real closed fields can also be characterized as follows: an $n$-ordered field $\calK$ is pseudo real closed if and only if every absolutely irreducible plane curve which has a simple point in
every real closure of $K$ has infinitely many $K$-rational points. Compare this characterization with the characterization of \emph{pseudo algebraically closed} fields which are those fields $M$ for which every absolutely irreducible plane curve has infinitely many $M$-rational points. The following theorem of Stone is essential in the study of pseudo real closed fields.
\begin{fact}[Stone]
\label{approxthm}
Let $\calK$ be an $n$-ordered field such that $<_i$ and $<_j$ induce different interval topologies for $1\leq i <j\leq n$ and let $I_i\subseteq K$ be an $<_i$-interval for $i=1,\ldots,n$. Then $\bigcap_{i=1}^nI_i \neq \emptyset$.
\end{fact}

\noindent We use $\overline{\calK}$ to denote the real closure of $K$ with respect to $<_1$ and we use \textbf{dense} to mean dense in the topology induced by $<_1$, unless otherwise specified. Let $\calLB = \{0,1,+,\cdot,-,<_1\}$, so $\calLB \subseteq \calLA$. Let $\TB$ be the $\calLB$-theory of real closed ordered fields,  let $\calL^2$, $T^2$ be as in Section \ref{s:setup}, and  let $T^d_n$ be the $\calL^2$-theory
$$T^2 \cup \big\{\forall y_1\forall y_2\exists x (A(x) \wedge y_1<_1 x <_1 y_2)\big\}.$$
The models of $T^d_n$ are real closed fields with a predicate for a dense pseudo real closed subfield with $n$ orderings, where the ordering on the bigger field agrees with the first ordering of the subfield. It is a fact that any model $\calK\models \TA$ is dense in $\overline{\calK}$, so the pair $(\overline{\calK},\calK)$ is a model of $T^d_n$. Also, if $\calR$ is a real closed ordered field containing $\overline{\calK}$ as a dense subfield, then $(\calR,\calK)$ is a model of $T^d_n$. The main result of this subsection is the following theorem:

\begin{thm}
\label{PRCML}
$T^d_n$ is an ML-theory.
\end{thm}

\noindent The proof of this theorem follows from the three lemmas below:

\begin{lem}
\label{PRCMLC}
$T^d_n$ satisfies the Mordell-Lang condition.
\begin{proof}
Let $(p,q,\phi,\psi)$ be a Mordell-Lang challenge, suppose that $\big((\calR,\calK),\vec{c},a\big)$ is a solution, and let $\big((\calR',\calK'), \vec{d}\ \big)$ be a contender. Since $\vec{c}\forkindep_{\vec{c}_\alpha}K$, we have that $a$ is in $\dcl(\vec{c}_\alpha)$. Thus, we may assume that $\phi(\vec{c}_\alpha,y)$ isolates the type $\tp_{\calLB}(a|\vec{c}_\alpha)$. Since $\big|\phi(\vec{c},\calK)\big| = 1$, it must be the case that $\phi(\vec{c}_\alpha,y)$ isolates the type $\tp_{\calLB}(a|\vec{c})$. We may also assume that $\phi$ is quantifier-free by quantifier elimination for real closed ordered fields. Since $\calLB\subseteq \calLA$ and since $\phi = \phi_A$, we have that $\phi(\vec{c}_\alpha,y)$ isolates the type $\tp_{\calLA}(a|\vec{c}_\alpha)$. Therefore $\big((\calR',\calK'),\vec{d},b\big)$ is also a solution for any $b \in \calK'$ such that $\calK' \models \phi(\vec{d}_\alpha,b)$.
\end{proof}
\end{lem}

\begin{lem}
\label{PRCnotsmall}
$T^d_n$ satisfies the codensity condition.
\begin{proof}
Let $(\calR,\calK) \models T^d_n$. By Lemma \ref{expODAG}, it suffices to show that $R$ is not $\calK$-small.
Let $P_1,\ldots,P_m \in R[X_1,\ldots X_k,X_{k+1}]$ be polynomials over $R$ and let
$$Z = \{z \in R:\text{ there are }a_1,\ldots,a_k \in K\text{ and }i \in \{1,\ldots,m\}\text{ such that }P_i(a_1,\ldots,a_k,z) = 0\}.$$
It suffices to show that $Z \neq R$. Let $\vec{c}$ be a tuple of elements in $R$ such that $P_1,\ldots,P_m \in K(\vec{c})[X_1,\ldots X_k,X_{k+1}]$ and let $d$ be the degree of the field extension $K(\vec{c})/K$. Let $e$ be the maximum degree of $X_{k+1}$ that appears in any of the $P_i$. Then the degree of $K(z)/K$ is at most $d+e$ for all $z \in Z$. We claim that $R$ contains elements of arbitrarily high degree over $K$, so $R$ cannot be equal to $Z$. Take $a \in K$ with $a>_1 0$ and $a<_2 0 $ (such an element exists by Fact \ref{approxthm}). Then for any $\ell = 1,2,\ldots$ there is $b \in R$ with $b^{2^\ell} = a$. An induction on $\ell$, using the fact that $a$ can't have any even roots in $K$, shows that $\deg\big(K(b)/K\big) = 2^\ell$ for such an element $b$.
\end{proof}
\end{lem}

\begin{lem}
\label{PRCdensity}
$T^d_n$ satisfies the density condition.
\begin{proof}
Fix a $\kappa$-saturated model $(\calR,\calK) \models T^d_n$ where $\kappa$ is uncountable. Fix $C \subseteq R$ with $|C|<\kappa$ and a non-algebraic unary $\calLB(C)$-type $q(x)$. By o-minimality of $\TB$, we may assume that $q$ is a cut in $\dcl(C)$. Let $p(x)$ be a unary $\calLA(C)$-type such that $q\models \qf(p|_{\calL})$. We show that $p_A \cup q$ is realizable (hence realized by saturation) in $(\calR,\calK)$. Consider the formula
$$\theta(x,\vec{c},b_1,b_2) := \big(A(x)\wedge \psi_A(\vec{c},x) \wedge( b_1 <_1x<_1b_2 ) \big)$$
where $\psi(\vec{c},x) \in p(x)$ (so $\vec{c}$ is a tuple from $A\cap C$) and $b_1$, $b_2 \in \dcl(C)$ with $q(x) \models b_1<_1x<_1b_2$. By Montenegro \cite[Theorem 3.13]{Mo17}, we can find quantifier-free $\calL(\vec{c})$-definable subsets $I_1,\ldots,I_\ell \subseteq K$ such that
\begin{itemize}
\item $I_k$ is $<_1$-convex and $<_1$-open for $k = 1,\ldots,\ell$,
\item $\psi(\vec{c},\calK)$ is dense in $I_k$ for $k=1\ldots,\ell$, and
\item $\psi(\vec{c},\calK) \setminus \big(\bigcup_{k=1}^\ell I_k\big)$ is a finite subset of $\dcl(\vec{c})$.
\end{itemize}
As $\qf(p|_{\calL})$ is non-algebraic, there is a unique $k \in \{1,\ldots,\ell\}$ such that $\qf(p|_{\calL})\models x \in I_k$. Now view $I_k$ as a subset of $R$ (defined by the same quantifier free $\calL(\vec{c})$-formula), so 
\[
\calR \models \exists x (x \in I_k\wedge b_1<_1x<_1b_2).
\]
 As $\psi_A(\vec{c},\calK)$ is dense in $I_k \cap K$,
it is also dense in $I_k$ and so 
\[
(\calR,\calK)\models \exists x \theta(x,\vec{c},b_1,b_2).\qedhere
\]
\end{proof}
\end{lem}

\noindent By \cite[Theorem 3.2.2]{vdD78}, the completions of $\TA$ are in bijective correspondence with the isomorphism classes of the fields of algebraic elements $\calK^{\operatorname{alg}}$ for models $\calK\models \TA$. Using this and Theorem \ref{PRCML}, we are able to characterize the completions of $T^d_n$:
\begin{cor}
Let $(\calR_1,\calK_1),(\calR_2,\calK_2) \models T^d_n$. The following are equivalent
\begin{enumerate}
\item $(\calR_1,\calK_1)\equiv(\calR_2,\calK_2)$,
\item $\calK_1\equiv \calK_2$,
\item $\calK_1^{\operatorname{alg}} \simeq \calK_2^{\operatorname{alg}}$.
\end{enumerate}
\end{cor}

\noindent Using the fact that $\TA$ is model-complete and that $\TB$ admits quantifier elimination, we have by Corollary \ref{nmc}:
\begin{cor}
$T^d_n$ is near model complete.
\end{cor}

\noindent The following corollary is immediate from Theorems \ref{opencore} and \ref{PRCML}.
\begin{cor}\label{thmc}
Every open set definable with parameters in a model of $T^d_n$ is semi-algebraic.
\end{cor}

\section{P-adics with a dense independent set}\label{s:padic}

In this section, let $p$ be prime and let $\TB$ be the theory of the $p$-adic field $\Q_p$ in the language $\calLB = \{0, 1, + \cdot, \calO, P_2,P_3,\ldots\}$ where $\calO$ is a unary predicate interpreted as the valuation ring of $\Q_p$ and $P_n$ is a unary predicate for every $n \geq 2$ with the interpretation $P_n(x) \Leftrightarrow \exists y (y^n = x)$. 
\begin{fact} The following fundamental facts about $T_{\beta}$ ensure the satisfaction of many of our conditions:
\begin{enumerate}
\item $\TB$ has quantifier elimination in the language $\calLB$ (due to Macintyre \cite{Ma76}).
\item Any infinite definable subset of a model of $\TB$ has nonempty interior with respect to the valuation topology (this follows from quantifier elimination).
\item The theory $\TB$ has definable Skolem functions. In particular, $\acl = \dcl$ in every model of $\TB$ (implicit in work of van den Dries \cite{vdD78}).
\end{enumerate}
\end{fact}

\noindent Let $\calLA '$ be a relational language disjoint from $\calLB$, and let $\TA '$ be a complete and consistent $\calLA '$-theory.
Let $\calLA$ be the expansion of $\calLA '$ by a binary predicate $E$ not already in $\calLB$ or $\calLA '$.
We now mirror the construction of pairs in \cite{HNW17}. For each $\calLA'$-formula $\phi$, we define an $\calLA$-formula $\phi_e$ as in \cite{HNW17}, that is, we replace every instance of equality ``$ x= y$'' in $\phi$ with ``$x E y$.''
We construct $\TA \supseteq \{ \theta_e : \theta \in \TA ' \}$ by requiring also that $E$ is an equivalence relation with infinite classes and that each relation $R$ in $\calLA$ is $E$-invariant. \newline

\noindent Let $T^2$ be as in Section \ref{s:setup} and let $T^* \supseteq T^2$ be the theory stating that in any model $(\calQ_p,\calA) \models T^*$:
\begin{itemize}
\item $A$ is dense in $Q_p$ with respect to the valuation topology and $\acl$-independent in $Q_p$,
\item Each equivalence class of $E$ is dense in $A$ with respect to the valuation topology.
\end{itemize}

\begin{lem}
The theory $T^*$ is consistent and $T^*$ interprets $\TA'$.
\begin{proof}
By the proof of \cite[1.11]{DMS16}, there exists a model $\calQ_p\models \TB$ and a family $(A_\gamma)_{\gamma < |\TA'|}$ of dense, pairwise disjoint $\acl$-independent subsets of $Q_p$ (one only needs to change ``open intervals'' to ``basic open balls''). By \cite[Lemma 2.2]{HNW17} (with $\TB$ in place of $T$), this model $\calQ_p$ admits an extension to a model $(\calQ_p,\calA) \models T^*$ (their proof of this lemma does not use o-minimality, so it goes through in our context). This shows that $T^*$ is consistent. To see that $T^*$ interprets $\TA'$, fix $(\calQ_p,\calA) \models T^*$, set $A' := A/E$, and expand $A'$ to a $\calLA'$-structure $\calA'$ by setting
\[
\calA' \models R\big([a_1]_E,\ldots,[a_n]_E\big):\Longleftrightarrow \calA \models R(a_1,\ldots,a_n)
\]
for each $a_1,\ldots,a_n \in A$ and each relation $R \in \calLA'$. This is well defined since every such $R$ is $E$-invarient. Of course, $\calA'$ is interpretable in $(\calQ_p,\calA)$ and since $ \TA \supseteq \{ \theta_e : \theta \in \TA ' \}$, we have that $\calA' \models \TA$. See \cite[Proposition 2.4]{HNW17} for additional details.
\end{proof}
\end{lem}

\begin{lem}
\label{p-adicnotsmall}
Let $(\calQ_p,\calA) \models T^*$. Then no open set in $\calQ_p$ is $\calA$-small.
\begin{proof}
It suffices to show that no basic open ball around 0 is $\calA$-small. Our argument is essentially \cite[2.1]{DMS16}. Let $v$ denote the valuation on $\calQ_p$ and let $\Gamma := v(Q_p^\times)$. Let $1_\Gamma:= v(p)$ be the least positive element of $\Gamma$. Suppose for contradiction that there is $r \in \Gamma$ such that the basic open ball $B_{r}:= \{x \in Q_p:v(x) > r\}$ is $\calA$-small.
Since $\TB$ has definable Skolem functions, there is a $\calLB(Q_p)$-definable function $g:Q_p^m \to Q_p$ such that $B_r \subseteq g(A^{m})$.
Take a finite set $D \subseteq Q_p$ such that $A \cup D$ is $\acl$-independent and such that $g$ is $\calLB(D)$-definable (this can be done by increasing $m$, since any defining parameters from $A$ can be viewed as variables).
Set $\ell= |D|+m+2$. Set $r_1:= r$ and set $r_i := r_{i-1}+1_\Gamma$ for each $1<i \leq \ell$.
By density, we can find elements $a_i \in A \cap ( B_{r_i}\setminus  B_{r_{i+1}})$ for each $1 \leq i < \ell$ and an element $a_\ell \in A\cap B_{r_\ell}$.
Set $d :=a_1+a_2 + \ldots+a_{\ell}$
and observe that
$$v (d) =v(a_1+a_2 + \ldots+a_{\ell}) \geq \min \big\{ v(a_1),\ldots,v(a_{\ell})\big \}>r_1,$$
so $d \in B_{r_1}$.
By our assumption we can write $d= g(\vec{c})$ for some tuple $\vec{c} \in A^{m}$.
Thus for each $i \in \{ 1, \ldots ,\ell \}$ we have that
\[
a_i \in \acl\big(D \cup\{c_1,\ldots,c_{m}\}\cup \{a_j:j \neq i\}\big).
\]
Since $A \cup D$ is $\acl$-independent, this means that $a_i \in D \cup\{c_1,\ldots,c_{m}\}\cup \{a_j:j \neq i\}$ for each $i \in \{ 1, \ldots ,\ell \}$. Since $\ell >  |D|+m+1$, it must be the case that $a_i = a_j$ for some $i<j\leq \ell$, contradicting our disjoint selections of $a_1,\ldots,a_\ell$.
\end{proof}
\end{lem}

\noindent We now appeal to the independence and topological density of the predicate subset to conclude the following:

\begin{thm}\label{padicML}
$T^*$ is an ML-theory.
\begin{proof}
To see that the density condition holds, we remark that for any $\kappa$-saturated model $(\calQ_p,\calA) \models T^*$ with $\kappa> |T^2|$, for any $C \subseteq Q_p$ with $|C|<\kappa$ and for any non-algebraic $\calLB(C)$-type $q(x)$, every formula in $q(x)$ defines a set with nonempty interior. Let $p$ be any $\calL(A\cap C)$-type such that $q \models \qf(p|_{\calL})$ and fix $a \in A$ realizing $p_A$. Fix also $\phi(x) \in q(x)$. By density of the equivalence classes of $E$ in $Q_p$, there is an element $a' \in A$ such that $a'$ is in the interior of the set defined by $\phi$ and such that $a' E a$ (thus $a'$ realizes $p_A$).
By saturation we may find an element in $A$ realizing both $p_A$ and $q$. 
The codensity condition follows from Lemma \ref{p-adicnotsmall}, the fact that every unary nonalgebraic $\calLB(C)$-formula defines a set with nonempty interior, and saturation. Since $T^*$ satisfies the density and codensity conditions and since $A$ is $\acl$-independent in every model of $T^*$, we have that $T^*$ is a particular example of the theory $T_{\ind}$. Thus, $T^*$ is an ML-theory by Proposition \ref{ind}.
\end{proof}
\end{thm}

\noindent We have the following consequences of Theorem \ref{padicML}, Corollaries \ref{complete} and \ref{nmc}, and Theorems \ref{opencore} and \ref{NIP}.
\begin{cor}
$T^*$ is complete. If $\TA'$ is model complete, then $T^*$ is near-model complete.
\end{cor}

\begin{cor}
Every open set definable with parameters in a model of $T^*$ is semi-algebraic.
\end{cor}

\begin{cor}
If $\TA'$ is NIP, then so is $T^*$.
\end{cor}


\end{document}